\documentclass[reqno,10pt]{amsart}
\usepackage {amsmath,amssymb}
\usepackage {color,graphicx}
\usepackage {times}
\usepackage {float}
\usepackage[dvips]{hyperref}

\newcommand {\PP}{{\mathbb P}}

\newcommand {\RR}{{\mathbb R}}
\newcommand {\ZZ}{{\mathbb Z}}

\newcommand {\calM}{{\mathcal M}}

\DeclareMathOperator{\ev}{ev}
\DeclareMathOperator{\val}{val}

\DeclareMathOperator{\Int}{Int}
\DeclareMathOperator{\codim}{codim}
\DeclareMathOperator{\ft}{ft}
\DeclareMathOperator{\divisor}{div}
\DeclareMathOperator{\mult}{mult}
\DeclareMathOperator{\rr}{rr}

\DeclareMathOperator {\Area}{Area}
\DeclareMathOperator {\lab}{lab}

\newcommand {\dunion}{\; \mbox {\raisebox{0.25ex}{$\cdot$} \kern-1.63ex $\cup$}  \,}

\newcommand {\twolines}[2]{\stackrel{\scriptstyle {#1}}{#2}}
\newcommand {\M}[2]{\calM_{0,#1}^{\text{lab}}(\RR^2,d)}
\newcommand {\dcup}{\hspace{-0.7ex}
                      \begin{array}{c}
                        \\[-3.12ex]
                        \cdot \\[-2.3ex]
                        \cup
                      \end{array}
                      \hspace{-0.7ex}
                   }

\newtheorem {theorem}{Theorem}[section]

\newtheorem {lemma}[theorem]{Lemma}

\newtheorem {propdef}[theorem]{Proposition and Definition}

\theoremstyle {definition}
\newtheorem {definition}[theorem]{Definition}
\newtheorem {example}[theorem]{Example}

\theoremstyle {remark}
\newtheorem {remark}[theorem]{Remark}
\newtheorem {construction}[theorem]{Construction}

\newtheorem {notation}[theorem]{Notation}

\begin {document}

\title [Tropical descendant Gromov-Witten invariants]{Tropical descendant Gromov-Witten invariants}

\author {Hannah Markwig}
\address {Hannah Markwig, CRCG, Mathematisches Institut, Georg-August-Universit\"at G\"ottingen, Bunsenstr. 3-5, 37073 G\"ottingen, Germany }
\email {hannah@uni-math.gwdg.de}

\author {Johannes Rau}
\address {Johannes Rau, Fachbereich Mathematik, Technische Universit\"at Kaiserslautern, Postfach 3049, 67653 Kaiserslautern, Germany}
\email {jrau@mathematik.uni-kl.de}

\subjclass[2000]{Primary 14N35, 52B20}

\begin{abstract}
We define tropical Psi-classes on $\calM_{0,n}(\RR^2,d)$ and 
consider intersection products of Psi-classes and pull-backs of evaluations on this space. We show a certain WDVV equation which is sufficient to prove that tropical numbers of curves satisfying certain Psi- and evaluation conditions are equal to the corresponding classical numbers.
We present an algorithm that generalizes Mikhalkin's lattice path algorithm and counts rational plane tropical curves satisfying certain Psi- and evaluation conditions.
\end{abstract}

\maketitle

\section {Introduction} \label {introduction}


Psi-classes $\psi_i$ are certain divisor classes on spaces of stable curves or stable maps, $\overline{M}_{g,n}$ or $\overline{M}_{g,n}(\PP^r,d)$, which arise as the first Chern class of the line bundle whose fiber over a point $(C,{x}_1,\ldots,x_n)$ (or $(C,x_1,\ldots,x_n,f)$) is the cotangent space of $C$ at $x_i$.
They are for example useful to count curves with tangency conditions.
To count curves that satisfy incidence conditions (e.g.\ pass through given points), one defines evaluation maps on the space of stable maps, $\ev_i:\overline{M}_{g,n}(\PP^r,d)\rightarrow \PP^r$ that send a stable map $(C,x_1,\ldots,x_n,f)$ to the image $f(x_i)$ of the marked point $i$. Then we can pullback the conditions via the evaluation map. Finally, we can intersect pullbacks of evaluation maps and Psi-classes on $\overline{M}_{g,n}(\PP^r,d)$. The degrees of such top-dimensional
intersection products are called descendant Gromov-Witten invariants and have been studied in Gromov-Witten theory. \\
The aim of this paper is to define tropical analogues of rational descendant Gromov-Witten invariants in the plane
and to show that, under certain assumptions on the distribution of the Psi- and evaluation conditions, they coincide with their conventional counterparts. \\
To do so, we use the constructions of moduli spaces of abstract and parameterized rational tropical curves
as tropical varieties and the definition of evaluation maps on the latter ones (\cite{SS04}, \cite{GKM07}, \cite{Mi07}). Moreover, \cite{Mi07} already defines Psi-classes on the space of abstract tropical curves $\mathcal{M}_{0,n}$
and \cite{KM07} deals with their intersections. In this paper, we define Psi-classes on
the space of parameterized tropical curves. Together with the intersection theory of \cite{Mi06} and \cite{AR07}, we have all tools to define descendant Gromov-Witten invariants. We show that these invariants are independent of the 
position and ``type'' of the evaluation conditions and that they fulfill the string and divisor equation.
Then we show that the invariants also fulfill a certain WDVV equation which can be
used to determine the numbers recursively. As the classical numbers fulfill the same equations, it
remains to compare the initial numbers appearing in the recursion to show that the classical and tropical
invariants coincide. \\
These results can only be achieved for invariants such that the Psi-conditions come together with point conditions, and not alone or with line conditions. Note that one should expect such restrictions
as we work with a non-compact moduli space that does \emph{not}
parametrize curves with components in the ``boundary'' of $\RR^n$. Hence the tropical descendant Gromov-Witten invariants are different from the classical ones in some cases, namely whenever tropical curves are ``missing'' in the corresponding tropical count. However, we show that this does not happen when Psi-conditions
always come together with point conditions. \\
To prove the WDVV equation we show that the weight of a curve in an intersection
product can be computed locally as the determinant of a map (which basically collects all evaluation maps)
and then proceed similarly to \cite{GM053}. Finally, we present a tropical algorithm similar to Mikhalkin's lattice path count (\cite{Mi03}) to determine the numbers of rational plane tropical curves passing through points and satisfying Psi-conditions.

Our definition of tropical descendant invariants partly agrees with Mark Gross' definition which was found independently in his study of mirror symmetry (\cite{G}). 

The paper is organized as follows. In section \ref{invariants}, we recall some facts about tropical moduli spaces and tropical intersection theory that we need. Then we define Psi-classes on the space of parametrized tropical curves and tropical descendant invariants. In section \ref{generalconditions} we define what it means for incidence conditions to be general and what consequences arise for our tropical descendant invariants if we choose the conditions to be general. In section \ref{contractededge}, we show that every tropical curve in an intersection product of Psi-classes, point and line evaluations \emph{and} the pullback of a point with a large coordinate in $\mathcal{M}_{0,4}$ under the forgetful map contains a contracted bounded edge. 
Thus the tropical curve can be interpreted as a reducible curve by cutting it along this contracted bounded edge. In section \ref{multiplicities} we show that the weight of a tropical curve in a zero-dimensional intersection product can be computed using a determinant of a linear matrix. We use this in section \ref{splittingcurves} to show that the weight of tropical curves with a contracted bounded edge can be (almost) split into two factors corresponding to the irreducible components. In section \ref{stringdivisor}, we show the string equation and the divisor equation for our tropical descendant invariants. In section \ref{recursion} finally, we collect our results to prove that our tropical descendant invariants satisfy a certain WDVV equation, and we conclude that the tropical invariants are equal to the corresponding classical invariants that satisfy the same recursion. In section \ref{latticepaths}, we describe an algorithm similar to the lattice path count that determines tropical descendant invariants.

We would like to thank A.~Gathmann and M.~Kerber for useful discussions.

\section {Defining the invariants} \label {invariants}


First of all, let us briefly recall the constructions from \cite{AR07}
that we need here: 

A \emph{cycle $X$} is a balanced (weighted, pure-dimensional, rational and 
polyhedral) complex (resp. fan) in $\RR^n$. The top-dimensional polyhedra (resp. cones)
in $X$ are called \emph{facets}, the codimension one polyhedra (resp. cones) are
called \emph{ridges}. The integer weights assigned to each facet $\sigma$ 
are denoted by $\omega(\sigma)$.
\emph{Balanced} means that the weighted
sum of the primitive vectors of the facets $\sigma$ around a ridge $\tau \in X$
$$
  \sum_{\twolines{\sigma \in X^{(\dim(X))}}{\tau < \sigma}} 
    \omega(\sigma) v_{\sigma / \tau}
$$
vanishes ``modulo $\tau$'', or, precisely, lies in the linear vector space
spanned by $\tau$, denoted by $V_\tau$.
Here, a \emph{primitive vector $v_{\sigma / \tau}$ of $\sigma$ modulo
$\tau$} is a integer vector in $\ZZ^n$ that points from $\tau$ towards
$\sigma$ and fulfills the primitive condition: The lattice
$\ZZ {v}_{\sigma / \tau} + (V_\tau \cap \ZZ^n)$ must be equal to the 
lattice $V_\sigma \cap \ZZ^n$. Slightly differently, in \cite{AR07}
the class of $v_{\sigma / \tau}$ modulo $V_\tau$ is called primitive
vector and $v_{\sigma / \tau}$ is just a representative of it. \\
For us, a polyhedron $\sigma$ is always understood to be closed. The 
\emph{(relative) interior} $\Int(\sigma)$ is the topological interior 
of $\sigma$ in
its affine span (e.g. $\Int(\{P\})= \{P\}$). The
\emph{support of $X$}, denoted by $|X|$, is the union of all facets
in $X$ with non-zero weight.

A \emph{(non-zero) rational function on
$X$} is a function $\varphi : |X| \rightarrow \RR$ 
that is affine (resp. linear) with rational 
slope on each polyhedron (resp. cone). 
The \emph{divisor of $\varphi$}, denoted by
$\divisor(\varphi) = \varphi \cdot X$, is the
balanced subcomplex (resp. subfan) of $X$ constructed in \cite[3.3]{AR07}, 
namely
the codimension one skeleton $X \setminus X^{(\dim X)}$ together with
the weights $\omega_\varphi(\tau)$ for each ridge
$\tau \in X$. 
These weights are given by the formula
$$
  \omega_\varphi(\tau) =
    \sum_{\twolines{\sigma \in X^{(\dim(X))}}{\tau < \sigma}} 
    \omega(\sigma) \varphi_\sigma(v_{\sigma / \tau}) 
    - \varphi_\tau\Big( \sum_{\twolines{\sigma \in X^{(\dim(X))}}{\tau < \sigma}} 
    \omega(\sigma) v_{\sigma / \tau}\Big),
$$
where ${\varphi}_\sigma : V_{\sigma} \rightarrow \RR$ denotes the 
linear part of the affine function
$\varphi|_\sigma$. Note that the balancing condition of $X$ around
$\tau$ ensures that the argument of $\varphi_\tau$ is an element
of $V_\tau$.
If $\varphi$ is globally affine (resp. linear), all weights
are zero, which we denote by $\varphi \cdot X = 0$. Let the \emph{support
of $\varphi$}, denoted by  $|\varphi|$, be the subcomplex of $X$ containing 
the points where $\varphi$ is not locally affine.
Then we have $|\varphi \cdot X| \subseteq |\varphi|$.
Furthermore, the
intersection product is bilinear (see \cite[3.6]{AR07}). 
As the restriction of a rational function to a subcycle is again a 
rational function, we can also form multiple intersection products
$\varphi_1 \cdot \ldots \cdot \varphi_l \cdot X$. In this case we will 
sometimes omit ``$\cdot X$'' to keep formulas shorter. Note that multiple
intersection products are commutative (see \cite[3.7]{AR07}).

By abuse of notation, a \emph{cycle} is also a class of balanced fans
with common refinement and agreeing weights. A \emph{rational function 
$\varphi$} on such a class
is just a rational function on a fan $X$ contained in the class.
We can generalize our
intersection product to such classes of fans $[X]$ by defining
$\varphi \cdot [X] := [\varphi \cdot X]$. In the following, we try
to avoid these technical aspects whenever possible. We will also
omit the brackets distinguishing between fans and their classes, hoping that
no confusion arises.

A \emph{morphism of cycles $X \subseteq \RR^n$ and $Y \subseteq \RR^m$}
is a map $f : |X| \rightarrow |Y|$ that comes from a linear
map from $\ZZ^n$ to $\ZZ^m$ and that maps each polyhedron (resp. cone)
of $X$ into one of $Y$.
Such a morphism \emph{pulls back rational functions} $\varphi$ on $Y$ 
to rational functions $f^*(\varphi) = \varphi \circ f$ on $X$. Note that
the second condition of a morphism, which is not required in \cite{AR07},
makes sure that we do not have to
refine $X$ further. $f^*(\varphi)$ is already affine (resp. linear)
on each cone. Furthermore, we can \emph{push forward subcycles} 
$Z$ of $X$ to subcycles $f_*(Z)$ of $Y$. This is due to
\cite[2.24 and 2.25]{GKM07} in the case of fans and can easily be generalized to complexes. We can omit further refinements here if we assume
that $f(\sigma) \in Y$ for all $\sigma \in X$. The \emph{projection
formula} (see \cite[4.8]{AR07}) connects all the above constructions via
$$
  f_*(f^*(\varphi) \cdot X) = \varphi \cdot f_*(X).
$$

Moreover, let us recall the basic facts of rational equivalence 
introduced in \cite[section 8]{AR07}. 
The \emph{degree} of a zero-dimensional cycle $Z$ is 
just the sum of all weights. Hence the push-forward of a
zero-dimensional cycle preserves degree. If $X$ is a one-dimensional cycle,
and $\varphi$ is a \emph{bounded} rational function, then
$\deg(\varphi \cdot X) = 0$ (see \cite[8.3]{AR07}). 
The pull-back of a bounded rational 
function is again bounded. Two functions
are called \emph{rationally equivalent} if they differ by the sum
of a bounded and a globally linear function. Hence (and by linearity 
of the pull-back) rational equivalence is preserved when pulled back.

An example for functions that are rationally equivalent
is given by translations
of functions on $\RR^n$.

\begin{lemma}[Translations are rationally equivalent] 
  \label{Translations}
  Let $h$ be a rational function on $\RR^n$, choose $v \in \RR^n$
  and consider $h'$ with $h'(x) := h(x + v)$. Then 
  $h$ and $h'$ are rationally equivalent.  
\end{lemma}

\begin{proof}
  Let $X$ be a subdivision of $\RR^n$ on which $h$ is a rational function.
  For each cone $\sigma \in X$, let $h_\sigma$ be the linear part of the affine
  function ${h|}_{\sigma}$. Take the maximum of the finitely many 
  $h_\sigma(v), \sigma \in X$ and call it $c$.
  Now, $X$ subdivides the line segment
  $x + \lambda v, \lambda \in [0,1]$ into $q$ line segments
  of length $\lambda_i$ contained in some polyhedron $\sigma_i$. This means
  $h(x + v)$ can be expressed as $h(x) + h_{\sigma_1}(\lambda_1 v) + \ldots +
  h_{\sigma_q}(\lambda_r v)$, where $\sum_i \lambda_i = 1$. This implies
  $$
    h(x + v) - h(x) \leq c,
  $$
  which proves that $h' - h$ is bounded.
\end{proof}

In the following, we will apply these constructions and results
to the case of Psi- and evaluation classes on the space of rational
plane curves.

The tropical analogue $\calM_{0,n}$ of the space of stable 
$n$-marked curves is the space of trees, or (a quotient of) 
the tropical Grassmanian (\cite{SS04}, \cite{GKM07}, \cite{Mi07}). 
Thus an abstract tropical curve is just a tree with $n$ marked 
ends and whose bounded edges $e$ are equipped with a length 
$l(e)\in \RR_{>0}$. The fan $\calM_{0,n}$ is stratified by 
cones corresponding to combinatorial types of trees. The facets 
correspond to $3$-valent trees. \\
The tropical analogue $\calM_{0,n}(\RR^2,d)$ of the space of 
stable maps has been studied in \cite{GKM07}.
An element of $\mathcal{M}_{0,n}(\RR^2,d)$ is an abstract 
tropical curve $\Gamma$ (i.e.\ a tree) together with a map 
$h:\Gamma\rightarrow \RR^2$ such that the image satisfies 
the balancing condition and marked ends are contracted to a 
point. An important feature of this definition 
is that it also allows to contract bounded edges,
as it will happen in section \ref{contractededge} and \ref{splittingcurves}.
If we furthermore also label the \emph{non-contracted} ends, we obtain the space $\M{n}{d}$.
The advantage of this space is that, after choosing the vertex of one marked 
end as root vertex, we can identify
$\M{n}{d}$ with $\calM_{n + 3d} \times \RR^2$, where the second 
factor describes the position of the root
vertex in $\RR^2$ (cf. \cite{GKM07}). 
In particular, in this sense $\M{n}{d}$ is a tropical variety. 
For enumerative purposes,
its difference to $\calM_{0,n}(\RR^2,d)$ cumulates in nothing 
but a factor $(d!)^3$ by which each invariant
in $\M{n}{d}$ must be divided to get the corresponding one in 
$\calM_{0,n}(\RR^2,d)$. 
Note that, independent 
of the choice of a root vertex, there exists a forgetful map 
$\ft' : \M{n}{d} \rightarrow \calM_{n + 3d}$ forgetting just the position
of the image of a curve in $\RR^2$. 
This forgetful map $\ft' : \M{n}{d} \rightarrow \calM_{n + 3d}$ is a morphism of
tropical varieties, as after choosing a root vertex and identifying $\M{n}{d}$ with
$\calM_{n + 3d} \times \RR^2$, $\ft'$ is just the projection onto the first factor.

Analogues of Psi-classes on tropical $\calM_{0,n}$ have been 
defined recently (\cite{Mi07}). $\psi_i$ with $i = 1, \ldots, n$ 
is the codimension one 
subcycle that consists of cones corresponding to trees where the 
marked end $i$ is at a $4$- or higher-valent vertex. How such Psi-classes 
intersect is discussed in \cite{KM07}. To do so, Psi-classes
$\psi_i, i=1,\ldots,n$ are defined as
divisors of rational functions ${f}_i$ on $\calM_{0,n}$ cf. 
\cite[proposition 3.5]{KM07}. As $\calM_{0,n}$ is simplicial, the function $f_i$ can be defined by 
specifying its values on the primitive vectors of the rays contained in $\calM_{0,n}$.
These rays are given by curves with only one bounded edge splitting up the marked ends into two
sets $I \dcup J = [n]$. Let $v_{I|J}$ be the corresponding primitive vector and assume w.l.o.g. 
$i \in I$, then $f_i$ is defined by
$$
  f_i(v_{I|J}) = \frac{|J| \, (|J|-1)}{(n-1)(n-2)}.
$$
Note that we denote by $f_i$ a multiple of what is called $f_i$ in \cite{KM07}, such that 
we obtain $\divisor(f_i)
= \psi_i$. We use these functions to pull back Psi-classes to $\M{n}{d}$.

\begin{definition}[Psi-classes for parameterized curves]
  For $i=1,\ldots,n$ we define \emph{the $i$-th Psi-class on $\M{n}{d}$}
  to be $\psi_i := \divisor(\ft'^*(f_i))$.
\end{definition}

\begin{remark}
  It can be shown that two rational functions on $\M{n}{d}$ (or 
  $\calM_{0,n}$) defining the same divisor cycle only differ by the
  restriction of a globally linear function. Hence, the choice
  of the functions defining our Psi-classes is not really important
  for intersection-theoretic purposes. This justifies that throughout 
  our paper we use the specific function $\ft^*(f_i)$
  to describe $\psi_i$ and in particular define
  $$
    \psi_i \cdot Y := \ft'^*(f_i) \cdot Y,
  $$
  where $Y$ is an arbitrary subcycle of $\M{n}{d}$. Note also that
  for our purposes we do not really need that the function describing
  $\psi_i$ is (nearly) unique. The only thing we need to know is
  contained in the following lemma.
\end{remark}

\begin{lemma}[Products of Psi-classes] \label{HowPsiProductsLookLike}
  Let $r_1, \ldots, r_n$ be positive integers and let
  $$X = \prod_{k = 1}^n \psi_k^{r_k} \cdot \M{n}{d}$$ be
  a product of Psi-classes. Then $X$ is the codimension-$\sum_k r_k$-subfan 
  of $\M{n}{d}$ consisting of cones $\sigma$
 corresponding to trees such that for each vertex $V$ we have $\val(V)=K(I_V)+3$, where $I_V$ denotes the set
$$I_V=\{k\in [n]\; : \text{ end }x_k \text{ is adjacent to } V\}\subset [n]$$
and $K(I)$ is a short notation for $K(I)=\sum_{k\in I}r_k$.
The weight of $\sigma$ equals
$$\omega(\sigma)=\frac{\prod_{V}K(I_V)!}{\prod_{k=1}^n r_k!}.$$
\end{lemma}

\begin{proof}
  Choose a root vertex and identify $\M{n}{d}$ with $\calM_{0, n+3d}
  \times \RR^2$. Then $\ft'$ is just the projection on the first
  factor and we can apply \cite[9.6]{AR07}, i.e. instead of 
  intersecting the pull-backs of the $f_k$ on the product, we can just
  intersect the $f_k$ on the first factor and then multiply with ${\RR}^2$.
  Thus,
  $$
    X = \left(\prod_{k = 1}^n \psi_k^{r_k} \cdot \calM_{0, n+3d}\right) \times \RR^2,
  $$
  where here $\psi_k$ denotes a Psi-class in $\calM_{0, n+3d}$.
  But in the case of non-parameterized curves, it is proved 
  in \cite[4.1]{KM07} that the valence of the vertices
  and the weights of the facets satisfy the formulas of the
  statement. Multiplying with $\RR^2$ does not disturb this,
  as the weight of $\RR^2$ is one and as the combinatorics
  of a curve remain unchanged under $\ft'$.
\end{proof}

\begin{remark}
  In particular the preceding lemma says that $\psi_i$ consists
  of those curves whose marked end $i$ is adjacent to an at least $4$-valent 
  vertex (where bounded edges as well as marked ends and non-contracted ends count towards the valence).
\end{remark}

Later on, we will also use the \emph{forgetful map}
$$
  \ft : \M{n}{d} \rightarrow \calM_{0,4},
$$
which forgets the map of a given curve $C$ to $\RR^2$ and all its ends but the first four marked ends (it also ``stabilizes'', which means
that, after forgetting one marked end, it replaces all two-valent vertices by straight edges while
adding up lengths). 

\begin{lemma}
  The forgetful map $\ft : \M{n}{d} \rightarrow \calM_{0,4}$ is a morphism of
  cycles.
\end{lemma}

\begin{proof}
	Let $\ft_n : \calM_{0,n} \rightarrow \calM_{0,n-1}$ be the forgetful
	map that just forgets the $i$-th end of an $i$-marked 
	non-parameterized curve. It is shown in \cite[3.9]{GKM07} that $\ft_n$ is
	a morphism for all integers $n \geq 4$. As mentioned above, the map
	$\ft'$ is a morphism, too. Thus, the statement follows from the formula
	$\ft = \ft_5 \circ \ldots \circ \ft_n \circ \ft'$.  
\end{proof}

Moreover, we use the \emph{evaluation maps} 
$$
  \ev_i : \M{n}{d} \rightarrow \RR^2
$$
assigning to a curve $C$ the position of its $i$-th marked end.
It is shown in \cite[4.8]{GKM07} that these maps are also morphisms
of cycles. Along these morphisms we will pull back lines and points.

\begin{definition}[Lines] \label{Lines}
  A \emph{line $G$} is a one-dimensional cycle in $\RR^2$ that is the divisor
  of a tropical polynomial of degree one. In other words, lines are 
  divisors of translations of the functions $\max\{x,y,0\}$, $\max\{x,0\}$, 
  $\max\{y,0\}$ or $\max\{x,y\}$.
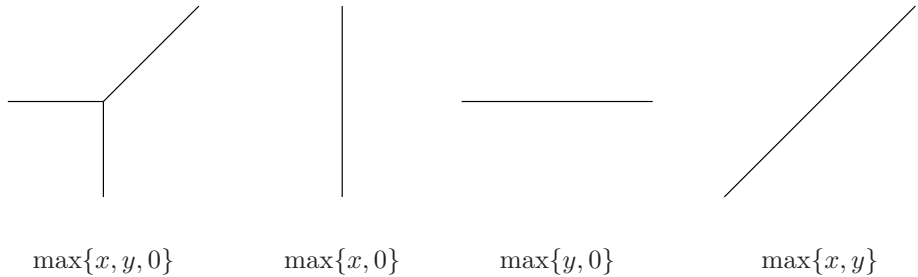
  \begin{center}
  \input{TypesOfLines.pstex_t}
  \end{center}
  Lines of type $\max\{x,y,0\}$ are also called \emph{non-degenerated}.
\end{definition}

We would like to pull back lines and points along an evaluation map
$\ev_i$. However, up to now, pull backs are only defined for functions,
not for cycles. Of course, we can choose rational functions cutting out
the line resp. point in question and pull them back instead. In the following
lemma we will show that,
for our purposes, the choice of describing functions plays no role.

\begin{notation}
We use the following notation: We have a total number of $l+m+n$
marked ends, which are subdivided into the three sets
$L \dunion M \dunion N = \{1, \ldots, l+m+n\}$, such that
$|L| = l$, $|M| = m$ and $|N| = n$. In the following, 
the ends $i \in L$ are unrestricted,
the ends $j \in M$ are restricted by lines $G_j$ (see \ref{Lines}) and 
the ends $k \in N$ have to meet points $P_k$. 
Furthermore we fix numbers $r_k, k \in N$ describing
how many Psi-classes we require at $k \in N$.
\end{notation}
\begin{lemma}
  Consider the intersection product
  $$
    Z:= \prod_{j \in M} \ev^*_j(G_j) 
    \prod_{k \in N} \ev^*_k(P_k) \psi_k^{r_k}
      \cdot \M{l+m+n}{d},
  $$
  where $\ev^*_j(G_j)$ stands for $\ev^*_j(h)$ with a 
  function $h$ cutting out $G_j$ and $\ev^*_k(P_k)$
  stands for $\ev^*_k(h_1) \cdot \ev^*_k(h_2)$ with
  function $h_1,h_2$ cutting out $P_k$.
  Then $Z$ is well-defined, i.e. it does not depend
  on the chosen rational functions. 
\end{lemma}

\begin{proof}
  Let $\ev := \ev_i$ be an 
  evaluation map and $G$ be a line. First we check
  that the intersection product $\ev^*(G) \cdot \M{l+m+n}{d}$ does not depend
  on the rational function describing $G$:
  Choose the vertex of the end $i$ as root vertex and identify $\M{l+m+n}{d}$ with
  $\calM_{0,l+m+n+3d} \times \RR^2$. Then $\ev$ is 
  just the projection onto the second factor. By
  \cite[9.6]{AR07} we deduce 
  $$
    \ev^*(G) \cdot (\calM_{0,l+m+n+3d} \times \RR^2)
      = \calM_{0,l+m+n+3d} \times G,
  $$
  which shows independence of the describing function. \\
  Now let $X = \varphi_1 \cdot \ldots \cdot \varphi_r \cdot \M{l+m+n}{d}$
  be a cycle given by arbitrary rational functions
  $\varphi_1, \ldots, \varphi_r$. Then, 
  by commutativity of the intersection product, the cycle
  $$
    \ev^*(G) \cdot X = 
      \varphi_1 \cdot \ldots \cdot \varphi_r \cdot \ev^*(G) \cdot \M{l+m+n}{d}
  $$
  is also well-defined. The same arguments work if we consider
  a point $P$ instead of $G$. But this suffices to conclude inductively
  that the big intersection product $Z$ is also well-defined.
  Moreover note that the same argument also shows that our choice
  of the function $f_i$ describing $\psi_i$ does not matter
  in this intersection product.
\end{proof}

We are now ready to define our tropical descendant Gromov-Witten invariants.

\begin{propdef} \label{InvariantsAndIndependence}
  Let $d,l,m,n$ and $r_k, k \in N$ be positive integers such that
  \begin{equation}\label{eq-dimension}
    l + m + n + 3d - 3 + 2 = m + 2n + \sum_{k \in N} r_k.
  \end{equation}
  Then we define the \emph{tropical descendant Gromov-Witten invariant} $\langle\tau_0(0)^l \tau_0(1)^m \prod_{k\in N} \tau_{r_k}(2)\rangle_d $ to be the number
  $$
   \langle\tau_0(0)^l \tau_0(1)^m \prod_{k\in N} \tau_{r_k}(2)\rangle_d :=
     \frac{1}{(d!)^3} \deg\left(\prod_{j \in M} \ev^*_j(G_j) 
    \prod_{k \in N} \ev^*_k(P_k) \psi_k^{r_k}
      \cdot \M{l+m+n}{d}\right).
  $$
  As indicated by the notation, this number only depends on $d,l,m,n,r_k, k \in N$,
  \emph{but not} on the lines $G_j$ and the points $P_k$.
\end{propdef}

\begin{proof}
  Lemma \ref{Translations} says
  that we can move around our points and lines arbitrarily,
  namely by translating the describing functions,
  without changing the degree. It remains
  to show that the type of the lines does not matter, for example the
  type of $G_1$. We will show that $\ev_1^*(G_1)\cdot F$ does not depend on the choice of the line $G_1$ for a one-dimensional cycle $F$, where 
  $$
    F =\prod_{j \in M \setminus \{1\}} \ev^*_j(G_j) 
    \prod_{k \in N} \ev^*_k(P_k) \psi_k^{r_k}
      \cdot \M{l+m+n}{d}
  $$
To see this, we have to use lemma \ref{PushForwardHasStandardDirections} 
which requires general conditions
and therefore is stated and proven in the next section of this article. 
It states that $\ev_{1*}(F)$ 
has only standard outer directions $-e_1, -e_2$ and
$e_1 + e_2$. Knowing this, we push forward $\ev_1^*(G_1)\cdot F$ via $\ev_1$, which does not 
change the degree and use the projection formula (\cite{AR07}). It 
tells us that $\ev_{1*}(\ev_1^*(G_1)\cdot F)=G_1\cdot \ev_{1*}(F)$. 
Now, as $\ev_{1*}(F)$ 
has only standard outer directions (at least for general conditions, which 
we can assume), any line intersects $\ev_{1*}(F)$ in the same number of points, 
not depending on the type. Note that 
lemma \ref{PushForwardHasStandardDirections} does not care about
the types of the lines appearing in the product of $F$. Thus we 
can apply the above argument inductively and see that
the types of all lines $G_j$ can be changed arbitrarily without 
changing the degree of the intersection product.
\end{proof}

\begin{remark}
The dimension of the space $\M{l+m+n}{d}=\calM_{0,m+n+l+3d}\times \RR^2$ is $l+m+n+3d-3+2$ since a $3$-valent tree with $m+l+n+3d$ ends has $l+m+n+3d-3$ bounded edges. The codimension of the intersection of Psi-classes is $\sum_{k \in N} r_k$. The pullback of a line has codimension $1$ and the pullback of a point codimension $2$. Hence the requirement (\ref{eq-dimension}) is equivalent to a $0$-dimensional expected dimension of the intersection.
\end{remark}

\begin{notation}
We will use the $\tau$-notation in a more general meaning:
  A product 
  \begin{displaymath}
    \big(\prod_{i \in I} \tau_{r_i}(c_i)\big)_d
  \end{displaymath} (with round brackets) stands for \emph{a
  cycle in $\M{|I|}{d}$}, obtained as the intersection product
  where we replace the $i$-th factor
  $\tau_{r_i}(c_i)$ by $\psi_i^{r_i} \ev^*(C_i)$. Here,
  $C_i$ is some point $P_i$ if $c_i = 2$, some line $G_i$ (of some type)
  if $c_i = 1$
  and the whole space $\RR^2$ (which means you can omit this 
  pull-back) if $c_i = 0$; thus $c_i$ describes the codimension
  of $C_i$. 
  If $\big(\prod_{i \in I} \tau_{r_i}(c_i)\big)_d$ is
  zero-dimensional, we denote, as before, by
\begin{displaymath}
  \langle\prod_{i \in I} \tau_{r_i}(c_i)\rangle_d =
     \frac{1}{(d!)^3} \deg\big(\prod_{i \in I} \tau_{r_i}(c_i)\big)_d
\end{displaymath}
the degree of the product above divided by $(d!)^3$. 
Note that a factor $\tau_0(0)$ can not be dropped
  in this notation as it stands for a marked end 
  that does not have to meet any condition at all. 
  \end{notation}
  
\begin{remark}\label{rem-lambdaind}
Later on, we will also allow the factor $\ft^*(\lambda)$
  in this notation, where $\lambda$ is an element in $\calM_{0,4}$
  and $\ft^*(\lambda)$ stands for the pull-back of a rational
  function on $\calM_{0,4}$ describing $\lambda$. Two such functions
  differ by an affine one, and so do the pull-backs. Hence, the
  intersection product containing $\ft^*(\lambda)$ as factor is still
  well-defined.
\end{remark}

\section {General incidence conditions} \label {generalconditions}


The invariants defined in \ref{InvariantsAndIndependence}
are well-defined also for ``special'' incidence conditions, e.g. if we choose
all points $P_i$ to coincide. In this case the \emph{set} of curves
fulfilling the conditions is of too big dimension,
but our intersection theory ensures
that the corresponding intersection product still has the correct dimension and degree.
However, many of the following arguments still require a notion
of ``general incidence conditions'' that ensures
that our intersection
product equals the set-theoretical count of curves 
fulfilling the incidence conditions (up to weights).

Let us start with the case of pulling back a single line in $\RR^2$.
Let $X$ be a subcomplex of $\M{n}{d}$, let 
$f : X \rightarrow \RR^2$ be a map that
is the restriction of a linear map (e.g.
morphisms like $f = \ev_i$) and let $G$ be a line in $\RR^2$.
Let $f^{-1}(G)$ be the subcomplex of $X$ containing all
polyhedra $\sigma \cap f^{-1}(\delta)$ for all $\sigma \in X$
and $\delta \in G$ (where $\delta$ denotes a cone in the polyhedral complex $G$).
Recall that the interior of a polyhedron $\Int(\sigma)$ denotes its topological
interior in its affine span.

\begin{lemma} \label{GeneralLine}
There exists a open dense subset $U\subset \RR^2$ such that for $v \in U$ and a translation
  $G' := G + v$ of $G$, it holds:
  \begin{enumerate}
    \item
      The subcomplex $f^{-1}(G')$ is either empty or of pure codimension $1$
      in $X$. 
    \item
      The interior of a facet of $f^{-1}(G')$
      is contained in the interior of a facet of $X$.
    \item
      For an element $C$ in the interior of a facet
      of $f^{-1}(G')$, the image $f(C)$ lies in the interior of a
      facet of $G$.
  \end{enumerate}
\end{lemma}

\begin{proof}
  Let $\sigma$ be a facet of $X$ and
  $\delta$ a one-dimensional polyhedron of $G$.
  Consider the map $q \circ f_\sigma$, where
  $f_\sigma : V_\sigma \rightarrow \RR^2$ is the extension
  of $f|_\sigma$ to $V_\sigma$ and $q : \RR^2 \rightarrow \RR^2/V_\delta$
  is the quotient map. This composition has either rank $1$ (in which
  case $\ker(q \circ f_\sigma)$ has codimension one in $V_\sigma$; hence,
  for a general translation $\delta'$ of $\delta$,
  the polyhedron $\sigma \cap f^{-1}(\delta')$ is either empty or 
  of codimension $1$ and intersecting the interior of $\sigma$) or
  has rank $0$ (then $\sigma \cap f^{-1}(\delta')$ is empty for 
  a general translation of $\delta'$ of $\delta$). As there are only
  finitely many pairs $\sigma, \delta$, the set of vectors $v \in \RR^2$
  such that these statements are true simultaneously is still open and dense. But note
  that all facets of $f^{-1}(G')$ can be obtained
  in this way for some pair $\sigma, \delta$. This shows
  part (a) and (b). \\
  Furthermore, let $V \in \RR^2$ be the vertex of $G$ (if $G$ is of type
  $\max\{x,y,0\}$). Applying the same argument to $V$ shows that for a general
  translation $V' := V + v$, the preimage $f^{-1}(V')$ has at least codimension
  $2$, which proves part $(c)$.
\end{proof}

\begin{definition} \label{GeneralConditions}
  Let $Z$ be an intersection product of the form 
  $(\prod_{i \in I} \tau_{r_i}(c_i))_d$ with incidence conditions
  $C_i$. Define $X:= \prod_{i \in I} \psi_i^{r_i} \cdot \M{|I|}{d}$.
  We call the conditions 
  \emph{general} if the following holds:
  \begin{enumerate}
		\item
			The subcomplex $S$ of $X$ containing all points $C \in X$
      fulfilling $\ev_i(C) \in C_i$ has dimension $\dim(S) = \dim(Z)$.
		\item
			The interior of a facet of $S$ 
			is contained in the interior of a facet of $X$.
		\item
			The interior of a facet $\sigma$ of $S$ 
			maps to the interior of a facet of $C_i$ under $\ev_i$. 
	  \item
	    Any intersection $C_i \cap C_j$, $i,j \in I$ has
	    expected codimension $c_i + c_j$.
  \end{enumerate}
\end{definition}

\begin{remark}\label{rem-sandz}
Let $S$ be the subcomplex of $X$ containing all the 
  curves $C \in X$ fulfilling $\ev_i(C) \in C_i$.
  Note that $Z$ is a subcomplex of $S$. This follows
  from the facts that the support of an intersection
  product is contained in the support of the intersecting
  rational function and that the support of a pull-back
  is contained in the preimage of the support of the
  pulled-back function. 
Note that in general we have $S=Z$ (as sets) if $\dim(S)=\dim(Z)$ is satisfied, the only thing that can happen in principle is that there are facets of $Z$ which get $0$ as a weight in the intersection product, although they are facets of $S$. For the intersection products we work with, this cannot happen though, since we only have a weight of $0$ if the set $S$ is of higher dimension (see section \ref{multiplicities}).
Hence for us the incidence conditions being general implies that $|Z|$ equals the set of curves satisfying the incidence conditions, and $\deg(Z)$ equals the number of curves satisfying the conditions, counted with weight.
\end{remark}

\begin{lemma}
  The set of general conditions in the space of all conditions (which can be
  identified with some big $\RR^N$ collecting all the translation vectors) is
  open and dense.
\end{lemma}

\begin{proof}
  The set of conditions fulfilling \ref{GeneralConditions} (d) is obviously open and dense.
  The remaining follows from recursively applying \ref{GeneralLine} to 
  $X$ and $\ev_1$, then $X \cap \ev_1^{-1}(C_1)$ and $\ev_2$, and so on.
  More precisely, if $c_i = 0$ we have nothing to do in this step, if 
  $C_i$ is a line, we apply \ref{GeneralLine}, and if $C_i$ is a point,
  we apply \ref{GeneralLine} twice for two lines intersecting set-theoretically
  in the single point $C_i$.  
\end{proof}

\begin{remark} \label{GeneralLambda}
  We also consider the following case:
  Let $X$ be a $1$-dimensional subcycle of $\M{n}{d}$ and consider
  the forgetful map $\ft : \M{n}{d} \rightarrow \calM_{0,4}$. We call 
  $\lambda \in \calM_{0,4}$ \emph{general}, if $\lambda \notin \ft(X^{(0)}) \cup
  \calM_{0,4}^{(0)}$, 
  where $X^{(0)}$ denotes the vertices
  of $X$ and $\calM_{0,4}^{(0)}$ denotes the single vertex of $\calM_{0,4}$. 
  This ensures that all points in $\ft|_X^{-1}(\lambda)$ lie in
  the interior of a one-dimensional polyhedron of $X$. 
\end{remark}

The following lemma describes the combinatorial type of the
curves which satisfy general incidence conditions.

\begin{lemma} \label{GeneralValence}
  Let $Z$ be an intersection product of the form
  $(\tau_0(0)^l \tau_0(1)^m \prod_{k\in N} \tau_{r_k}(2))_d$
  with general conditions. Then 
  \begin{enumerate}
		\item[(b')]
			For a curve $C$ in the interior of a facet 
			the following holds: All ends $k \in M \cup N$ 
			lie at different vertices and the valence of a vertex
			is $r_k + 3$ if $k \in N$ is adjacent to it and $3$ otherwise.
	\end{enumerate}
\end{lemma}

\begin{proof}
  Because of remark \ref{rem-sandz} we know that $Z\subset S$. In addition, condition \ref{GeneralConditions}
  (a) says that $Z$ and $S$ have the same dimension and
  therefore (b) and (c) also hold for curves in the interior
  of a facet of $Z$. \\
  Let $C$ be in the interior of a facet of $Z$.
  Condition \ref{GeneralConditions} (d) implies that
  $\ev_i(C) \neq \ev_j(C)$ for all $i \in M \cup N$,
  $j \in N$, as in this case $C_i \cap P_j$ is empty.
  If $i,j \in M$ would lie at the same vertex this
  would induce either a contracted bounded edge
  (which contradicts \ref{GeneralConditions} (a))
  or valence greater than $3$ of this vertex 
  (which contradicts \ref{GeneralConditions} (b)).
  Hence all ends in $M \cup N$ must lie at different
  vertices. The statement about the valence of
  the vertices follows from \ref{GeneralConditions} (b) and
  the description of $X$ in \ref{HowPsiProductsLookLike}.
\end{proof}

As a first application of our notion of general conditions we can now prove
the lemma which we promised and needed in the independence statement
\ref{InvariantsAndIndependence}. 

\begin{lemma} \label{PushForwardHasStandardDirections}
 Let $F$ be a one-dimensional cycle of the form
 $(\tau_0(0) \tau_0(0)^l \tau_0(1)^m \prod_{k\in N} \tau_{r_k}(2))_d$
 with general conditions. Let $x$ denote
 the marked end corresponding to the first factor $\tau_0(0)$.
 Then all of the unbounded rays of the push-forward
 $\ev_{x*}(F)$ have standard directions $-e_1, -e_2$ and
 $e_1 + e_2$.
\end{lemma}

\begin{proof}
 Let $\sigma$ be a facet of $F$. For a curve in the interior of $\sigma$ two possibilities can occur:
 Either $x$ is adjacent to a higher-valent vertex. Then by \ref{GeneralValence}
 also an end $k \in N$ interpolating the point $P_k$ lies at this vertex. Therefore,
 $\ev_x(\sigma) = \ev_k(\sigma) = \{P_k\}$. \\
 Secondly, $x$ might be adjacent to a $3$-valent vertex. 
Since $x$ itself is contracted, the two other edges which are adjacent are mapped to lines with opposite direction (because of the balancing condition). That means locally the image looks like a straight line with the marked point $h(x)$ on it. We can deform a curve in $\sigma$ in a one-dimensional family (thus covering $\sigma$) by changing the length of the two adjacent edges and thus making the point $h(x)$ move on the line. This movement is unbounded if and only if
 one of these two edges
 is an end. But then $\ev_x(\sigma)$ points to the same direction as this
 end, which is by definition one of the standard directions or $0$.
\end{proof}

\section {Contracted edges} \label {contractededge}


Let $F$ be a one-dimensional cycle of the form
$(\tau_0(0)^l \tau_0(1)^m \prod_{k\in N} \tau_{r_k}(2))_d$
with general conditions.
Remember that this implies that $|F|$ equals the set of curves satisfying the conditions.
\begin{notation}
We fix the type of the first four ends in the sense that we assume from now on
$1 \in L$, $2 \in M$ and $3,4 \in N$. 
\end{notation}
As before we denote by $\ft$ the forgetful map 
$\ft : \calM_{0,l+m+n}^{\text{lab}}(\RR^2, d) \rightarrow \calM_{0,4}$, 
which forgets the embedding and all ends but the first four 
marked ends.
It is the aim of this section to show
that for a very large $\calM_{0,4}$-coordinate $\lambda$, 
the curves in $\ft^{-1}(\lambda)\cap F$ (i.e. curves with such a large $\calM_{0,4}$-coordinate)
must contain a contracted bounded edge. We will use the contracted bounded edge in section \ref{splittingcurves} to split such curves into two components.

\begin{definition}
  Let $C$ be a curve in $\M{l+m+n}{d}$. For two 
	different marked ends $i_1,i_2$, we denote by 
	$S(i_1,i_2)$ the smallest connected subgraph of $C$ 
	containing $i_1$ and $i_2$ and call it \emph{the string of 
	$i_1$ and $i_2$}. Such a string $S(i_1,i_2)$ is called \emph{movable} 
	if $i_1,i_2 \in L \cup E$, where $E$ denotes the set of 
	non-contracted ends, and if $S(i_1,i_2)$ does not intersect 
	(the closure) of any $k$ for $k \in N$.
\end{definition}

\begin{lemma} \label{lem-string}
	Let $C$ be a curve in the interior of a facet of $F$. Then
	$C$ contains a movable string $S$.
\end{lemma}

\begin{proof}
	We know $\dim(F) = 1$, $\codim(F) = m + 2n + \sum_{k \in N} r_k$ 
	and $\dim(\M{l+m+n}{d}) = l + m + n + 3d  - 3 + 2$. Plugging in 
	all this in $\dim(F) + \codim(F) = \dim(\M{l+m+n}{d})$ leads to
	$$
		l + 3d = n + \sum_{k \in N} r_k + 2.
	$$
	On the other hand we can compute the number of connected 
	components of $\Gamma \setminus \bigcup_{k \in N} \bar{k}$: 
	Removing $k$ increases the number of connected components 
	by $r_k + 1$ as the valence of the adjacent vertex is 
	$r_k + 3$ by \ref{GeneralValence}. So, after removing all 
	$n$ ends, we arrive at $1 + n + \sum_{k \in N} r_k$ connected 
	components. The above equation tells us that 
	there is one more end in $L \cup E$ then there are
	connected components and therefore at least two ends $i_1, i_2 \in L \cup E$ 
	lie in the same component. Hence $S(i_1,i_2)$ 
	is a movable string.
\end{proof}

	By construction all vertices of a movable string are 3-valent.

\begin{lemma} \label{ContractedEdge}
	Let $\sigma$ be a facet of $F$ such that the corresponding interior curves do \emph{not}
	contain a contracted bounded edge. Then the image of $\sigma$ under $\ft$ is bounded.
\end{lemma}

\begin{proof}
	Let $C$ be a curve in the interior of $\sigma$. 
We will deform $C$ in a one-dimensional family inside $\sigma$. Since $\sigma$ is one-dimensional itself, this family covers $\sigma$.
By lemma \ref{lem-string} there exists a 
	movable string $S$ in $C$. 
  In the following, we show that 
	either $\sigma$ is bounded (i.e. the deformation of $C$ is bounded) or $\ft$ is constant on $\sigma$ (i.e. the deformation of $C$ does not affect $\ft$). \\
	Let $V$ be a vertex in $S$. We call $V$ \emph{degenerated} if we can deform $C$ one-dimensionally locally around $V$, i.e.\ if
	\begin{enumerate}
	  \item
	  either one of the adjacent edges is a marked end $i\in  L$,
	  \item
	   or one of the adjacent edges is a marked end $j \in M$
	    and the linear spans of the corresponding line $G_j$ at
	    $\ev_j(C)$ and of the other two edges adjacent to $V$
	    coincide (i.e. if the curve $C$ and the line $G$ do not
	    intersect transversally at $\ev_j(C)$),
	  \item
	    or all edges adjacent to $V$ are non-contracted, but their
	    span near $V$ is still only one-dimensional; w.l.o.g. we denote
	    the edge alone on one side of $V$ by $v$ and the two edges on
	    the other side by $v_1, v_2$.
  \end{enumerate}  
  \begin{center}
    \input{DegeneratedVertex.pstex_t}
  \end{center}
  If such a degenerated vertex exists, the $1$-dimensional deformation of the 
  curves inside $\sigma$ is given by moving this vertex and changing the lengths
  of the adjacent edges accordingly. We show that this
  movement is either bounded or, if not, the changed lengths do not influence $\ft$. \\
  Consider the cases (a) and (b) and let $v_1, v_2$ be the two other
  edges adjacent to $V$. At least one of the two edges, say $v_1$, is bounded.
  Then the movement is unbounded only if $v_2$ is unbounded. But $v_2$ cannot be
  contracted, as then $v_1$ would also be contracted (\emph{and} bounded). But this means
  that $\ft$ forgets $v_2$ and therefore also the length of $v_1$. \\
  Now consider the case (c). The balancing condition says $v = v_1 + v_2$ (by abuse
  of notation we denote the direction vectors by the same letters as the edges), which in
  particular implies that $v$ is not primitive and hence the edge $v$ has to be bounded. Now again, if
  we require the movement of $V$ to be unbounded, $v_1$ and $v_2$ must be unbounded.
  But they are also non-contracted which means that $\ft$ forgets them and the length 
  of $v$. \\
  So we are left with the case that all vertices of $S$ are non-degenerated.
 We can still describe the deformation of the
	curves inside $\sigma$ using the movement of the string: Take one of the ends
	of the string (which is necessarily non-contracted) and move it slightly
	in a non-zero direction modulo its linear span. Consider the next vertex $V$ and
	let $v$ be the adjacent edge not contained in the string. Then two things can
	happen: \\
	If $v$ is non-contracted (case A), our moved end will meet the affine span of
	$v$ at some point $P$ (as $V$ is non-degenerated). So we change the length of $v$ such that
	it ends at $P$ (while keeping the position of its second vertex fixed). Then we also move
	the second edge of the string to $P$ and go on to the next vertex. \\
	If $v$ is contracted (case B), our assumptions ensure that it is a marked end $j \in M$ and that
	the corresponding line $G_j$ intersects our curve transversally at $V$. Thus our moved
	edge will again meet $G_j$ at some point and by changing the
	lengths of the adjacent edges appropriately, the obtained curve will still meet $G_j$. \\
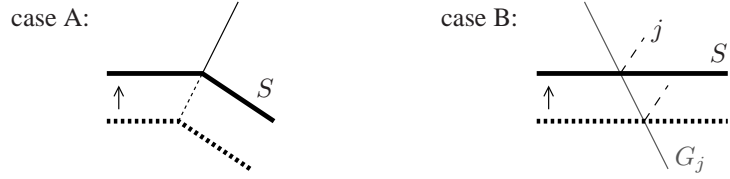
  \begin{center}
	  \input{MovingString.pstex_t}
	\end{center}
  In this way we can make our way through the string and finally obtain a deformation
	of the whole curve. Note that the non-degeneracy of all the vertices ensures that all edges
	of the string \emph{must} change their positions modulo their linear span and, hence, that all edges 
	adjacent to, but not contained in the string \emph{must} change their length. 
	In particular this means we cannot have more non-contracted
	ends adjacent to our string: Then we would have two different 
	strings providing two independent deformations
	of the curves inside $\sigma$, which is a contradiction 
	as $\sigma$ is 
	one-dimensional. \\
	Let us summarize: Our string $S$ is generated by two unique non-contracted
	ends $i_1, i_2$, all of its vertices are $3$-valent and the adjacent edges
	not contained in the string are either bounded or marked ends in $M$, where
	the corresponding line $G_j$ intersects transversally. The deformation only
	moves the string $S$; the adjacent edges are shortened or elongated and
	the other parts of the curve remain fixed. We want to show that, even if the 
	movement is unbounded, the considered $\calM_{0,4}$-coordinate is bounded. \\
	If
  there are bounded edges adjacent to $S$ to both sides of $S$ as in
  picture (a) below then the movement of the string is bounded. (This is true because if we move the string to either side, we can only move until the length of one of the adjacent bounded edges shrinks to $0$.) 
So we only have to consider the case when all adjacent bounded edges of $ S $
  are on the same side of $S $, say on the right side as in picture (b) below. Label the edges of $S$ (respectively, their direction vectors) by $ v_1,\dots,v_k $ and the adjacent bounded
  edges of the curve by $ w_1,\dots,w_{k-1} $ as in the picture. As above the
  movement of the string to the right is bounded. If
  one of the directions $ w_{i+1} $ is obtained from $ w_i $ by a left turn (as
  it is the case for $ i=1 $ in the picture) then the edges $ w_i $ and $
  w_{i+1} $ meet on the left of $S$. This restricts the movement of the string to the left, too, since the corresponding
  edge $ v_{i+1} $ then shrinks to length $0$. 
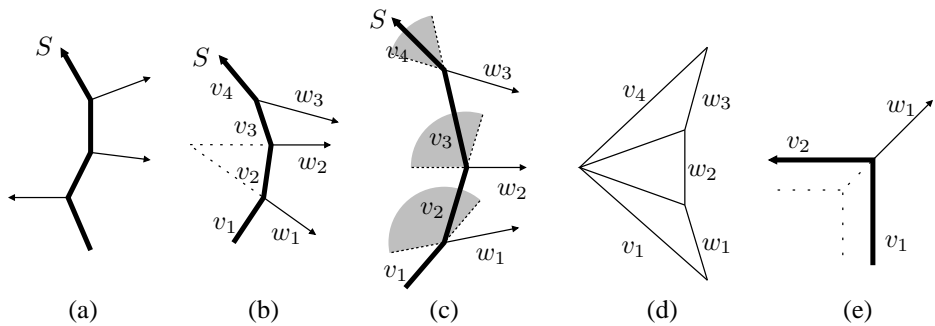
\begin{center}
\input{sides.pstex_t}
\end{center}
  So we can assume that for all $i$ the direction $ w_{i+1} $ is
  either the same as $ w_i $ or obtained from $ w_i $ by a right turn as in
  picture (c). The balancing condition then shows that for all $i$ both the
  directions $ v_{i+1} $ and $ -w_{i+1} $ lie in the angle between $ v_i $ and
  $ -w_i $ (shaded in the picture above). Therefore, all directions $ v_i
  $ and $ -w_i $ lie within the angle between $ v_1 $ and $ -w_1 $. In
  particular, the image of the string $ S $ cannot have any self-intersections in $
  \RR^2 $. We can therefore pass to the (local) dual picture (d)  where the edges dual to $ w_i $ correspond to a
  concave side of the polygon whose other two edges are the ones dual to $ v_1
  $ and $ v_k $.

  But note that there are no such concave polygons \emph {with integer
  vertices} if the two outer edges (dual to $ v_1 $ and $ v_k $) are two of the
  vectors $ \pm (1,0) $, $ \pm (0,1) $, $ \pm (1,-1) $ that can occur as dual
  edges of an end of a plane tropical curve of degree $d$. Therefore the string can consist at most of the two
  ends $ i_1 $ and $ i_2 $ that are connected to the rest of the curve
  by exactly one bounded edge $ w_1 $. This situation is shown in picture
  (e).

  In this case the movement of the string is indeed not bounded to the left.
  Note that then $ w_1 $ is the only internal edge whose length is not
  bounded. But by our assumptions $1$, $3$ and $4$ cannot lie 
  on $S$ but must lie on the other side of $w_1$; hence its
  length does not influence $\ft$. This finishes the proof.

\end{proof}

\section {Computing weights} \label{multiplicities}


In this section we prove that the weight of a curve
in a zero-dimensional intersection product can be computed
as the (absolute value of a) determinant of the linear map
that basically collects all the evaluation morphisms (and
forgetful morphism if present).
We will use this to express the weight of a reducible
curve (in the sense that it contains a contracted bounded edge)
in terms of the weights of the two components.

The following statement connects intersection products
with determinants.

Let $V = \RR \otimes \Lambda$ be a real vector space of dimension
$n$ with underlying lattice $\Lambda$ and let 
$h_1, \ldots, h_n \in \Lambda^\vee$
be linear functions on $\Lambda$ resp. $V$. By $H : V \rightarrow \RR^n$
we denote the linear map given by $x \mapsto (h_1(x), \ldots, h_n(x))$.
Choose lattice bases of $\Lambda$ and $\ZZ^n$ and consider the matrix
representation of $H$ with respect to these bases. Obviously, the absolute
value of the determinant of this matrix is independent of the choice of
bases; hence we denote it by $|\det(H)|$. \\
On the other hand, we can consider the rational functions 
$\varphi_i = \max\{h_i,0\}$ on $V$. To do so, we give $V$ the
fan structure consisting of all cones on which each $h_i$ is either
positive or zero or negative. These rational functions form a 
zero-dimensional intersection product, which obviously consists of only
$\{0\}$ with a certain weight.

\begin{lemma} \label{MultEqualDetLocal}
  The weight of $\{0\}$ appearing in $\varphi_n \cdot \ldots \cdot
  \varphi_1 \cdot V$ is equal to $|\det(H)|$.
\end{lemma}

\begin{proof}
  Let us first assume that ${h}_1, \ldots, h_n$ form a lattice basis
  of $\Lambda^\vee$ with dual basis $\tilde{h}_1, \ldots, \tilde{h}_n$. \\
  We compute $\varphi_1
  \cdot V$: For each ridge $\tau$ of $V$ (with the fan structure described
  above) there exists a unique $j$ such that $\tau \subseteq h_j^\bot$.
  Then there are two facets containing $\tau$
  (where $h_j$ is positive resp. negative) 
  and the corresponding (representatives of the) primitive vectors
  are $\tilde{h}_j$ resp. $-\tilde{h}_j$. Therefore the weight of $\tau$ in the 
  intersection product $\varphi_1 \cdot V$ can be computed as
  $$
    \omega_{\varphi_1}(\tau) =
      \max\{h_1(\tilde{h}_j),0\} + \max\{h_1(-\tilde{h}_j),0\} =
      h_1(\tilde{h}_j) + 0 = h_1(\tilde{h}_j).
  $$
  Hence, when omitting cones with weight $0$, $\varphi_1 \cdot V$ consists
  of all cones contained in $h_1^\bot$ and all weights are $1$. Now we can 
  apply induction on $h_2|_{h_1^\bot},  \ldots, h_n|_{h_1^\bot}$ and conclude
  that $\varphi_n \cdot \ldots \cdot \varphi_1 \cdot V$ produces $\{0\}$ with
  weight $1$. \\
  On the other hand, the matrix representation of $H$ with respect to the
  basis $\tilde{h}_1, \ldots, \tilde{h}_n$ for $\Lambda$ and the standard basis for $\ZZ^n$
  is just the unit matrix. Hence, $|\det(H)| = 1$, which proves the statement
  in the special case. \\
  General case: For general $h_1, \ldots, h_n$ we can choose a lattice basis
  $l_1, \ldots, l_n$ of $\Lambda^\vee$ such that
  $$
    \begin{array}{rcl}
      h_1 & = & a_{1,1} l_1, \\
      h_2 & = & a_{2,1} l_1 + a_{2,2} l_2, \\
      & \vdots & \\
      h_n & = & a_{n,1} l_1 + \ldots + a_{n,n} l_n,
    \end{array}
  $$
  where the $a_{i,j}$ are integers. Then we get
  $|\det(H)| = |\det((a_{i,j}))| = |a_{1,1} \cdot \ldots \cdot a_{n,n}|$. \\
  On the other hand, let us compute 
  that $\varphi_n \cdot \ldots \cdot \varphi_1 \cdot V$
  produces $\{0\}$ with weight $|a_{1,1} \cdot \ldots \cdot a_{n,n}|$.
  We saw in the special case that $\max\{l_1,0\} \cdot V$ is $l_1^\bot$
  with weight $1$. The above equations tell us that $\max\{h_1,0\} = 
  |a_{1,1}| \cdot \max\{l_1,0\}$. Using the linearity of the 
  intersection product,
  we deduce that $\varphi_1 \cdot V$ is $l_1^\bot$ with weight $|a_{1,1}|$.
  Now we apply induction on $n$ again: After restricting all functions to $l_1^\bot$, 
  we can omit all terms $a_{i,1} l_1$ in the above equations (in particular
  we can omit the first equation). Hence we can apply our induction
  hypothesis and conclude that
  $\varphi_n \cdot \ldots \cdot \varphi_2 \cdot l_1^\bot$ produces 
  $\{0\}$ with weight $|a_{2,2} \cdot \ldots \cdot a_{n,n}|$. Hence,
  $\varphi_n \cdot \ldots \cdot \varphi_2 \cdot (\varphi_1 \cdot V) =
  |a_{1,1}| \cdot \varphi_n \cdot \ldots \cdot \varphi_2 \cdot l_1^\bot
  = \{0\}$ with weight $|a_{1,1} \cdot \ldots \cdot a_{n,n}|$.
\end{proof}

With this tool we can express all weights occurring in a
zero-dimensional intersection product in terms of absolute values
of determinants.

\begin{notation}\label{not-z}
Let $Z$ be a zero-dimensional intersection product of the 
form $(\tau_0(0)^l \tau_0(1)^m$ \linebreak $\prod_{k\in N} \tau_{r_k}(2))_d$ or
of the form $(\ft^*(\lambda) \cdot 
\tau_0(0)^l \tau_0(1)^m \prod_{k\in N} \tau_{r_k}(2))_d$
with general condition $G_j$ and $P_k$ (and $\lambda$, resp.). The
set of curves $S$ that fulfill the incidence conditions set-theoretically
is finite by \ref{GeneralConditions} (a). 
For simplicity, let us furthermore assume that all points $P_k = (p_1,p_2)$ 
are described by the rational functions 
$\max\{x,p_1\}$ and $\max\{y,p_2\}$ on $\RR^2$
and that all lines $G_j$ are vertical, i.e.\ of type $\max\{x,0\}$ (i.e. are
given by a rational function $\max\{x,c_j\}$). 
$\lambda$ can
be described by $\max\{x,\lambda\}$, where $x$ is the coordinate of the ray
in whose interior $\lambda$ lies (see \ref{GeneralLambda}).

Denote $X := \prod_{k \in N} \psi_k^{r_k} \cdot \M{n}{d}$. 
We then consider the morphisms
\begin{eqnarray*}
  \ev : X & \rightarrow & \RR^m \times (\RR^2)^n, \\
        C & \mapsto     & 
          \Big(\big(\ev_j (C)_x\big)_{j \in M}, 
          \big(\ev_k (C)\big)_{k \in N}\Big),
\end{eqnarray*}
respectively
\begin{eqnarray*}
  \ft \times \ev : X & \rightarrow & \calM_{0,4} \times \RR^m \times (\RR^2)^n, \\
        C & \mapsto     & 
          \Big((\ft(C), \big(\ev_j (C)_x\big)_{j \in M}, 
          \big(\ev_k (C)\big)_{k \in N}\Big),
\end{eqnarray*}
where $\ev_j (C)_x$ denotes the first coordinate of the point
$\ev_j C \in \RR^2$. Thus, these morphisms evaluate at each end 
$i \in M \dunion N$
and keep all coordinates if $i \in N$ and 
only the first coordinate if $i \in M$. \\
Let $C$ be a curve in the interior of a facet
$\sigma$ of $X$ (and with $\ft(C)$ not being the vertex of $\calM_{0,4}$). 
Then $\ev$ (resp. $\ft \times \ev$) is affine in a neighborhood of $C$ and 
we define $|{\det}_{C}(\ev)|$ (resp. $|\det_C(\ft \times \ev)|$) to be
$|\det(H)|$, where $H$ is the linear part of $\ev$ (resp. $\ft \times
\ev$) at $C$.
\end{notation}
\begin{theorem} \label{MultEqualDet}
  The zero-dimensional intersection product $Z$ (as in notation \ref{not-z}) can be computed as
  $$
    Z = \sum_{C \in S} |{\det}_C(\ev)| \cdot C,
  $$
  resp. 
  $$
    Z = \sum_{C \in S} |{\det}_C(\ft \times \ev)| \cdot C,
  $$  
  i.e. the weight of a curve $C \in Z$ is just $|{\det}_C(\ev)|$
  (resp. $|\det_C(\ft \times \ev)|$).
\end{theorem}

\begin{proof}
  Each $C \in Z$ is contained in the interior
  of a facet $\sigma$ of $X$ (see \ref{GeneralConditions} (a)).
  In \ref{HowPsiProductsLookLike} the weight of $\sigma$ in $X$ 
  was computed to be 
  $$
    \omega(\sigma)=\frac{\prod_{V}K(I_V)!}{\prod_{k=1}^n r_k!}.
  $$
  But we know from \ref{GeneralValence} that no two marked ends
  lie at a common vertex and hence 
  $$
    K(I_V) = \left\{\begin{array}{cl}
                    	r_k & \text{if $k$ is adjacent to $V$,} \\
                    	0   & \text{otherwise.}
                    \end{array}
             \right.
  $$                  
  Therefore we can cancel the fraction defining $\omega(\sigma)$ down to $1$. \\
  As the computation of the weight of $C$ is local, we can replace
  $X$ by $V := \RR \cdot \sigma$. On the other hand, locally around $C$,
  all the pull backs along $\ev$ and $\ft$ are of the form $\max\{a,c\}$, where
  $a$ is an affine function on $V$ and $c$ is a constant.
  To be more precise, $a$ is exactly one of the coordinate functions
  of $\ev$ (resp. $\ft \times \ev$). 
  Now, up to translations and subtracting constant terms,
  we are in the situation of \ref{MultEqualDetLocal}:
  The weight of $C$ equals the absolute value of the determinant of
  linear map $H$ whose coordinate functions are
  the linear parts of the affine functions $a$. Thus $H$
  is the linear part of $\ev$ (resp. $\ft \times \ev$) on $V$, and we can
  conclude that the weight of $C$ in $Z$ is precisely $|\det(H)| =|\det_C(\ev)|$
  (resp. $= |\det_C(\ft \times \ev)|$).
\end{proof}

\begin{remark}
  If we dropped the requirement that Psi-conditions
  are only allowed at marked ends that are also restricted
  by a point condition, we could still prove a formula similar to
  the above one:
  For a zero-dimensional intersection product of arbitrary Psi- and
  evaluation classes, the weights can still be computed as the absolute
  value of an appropriate determinant times the weight of
  the corresponding facet in $X$. In particular, this shows that
  such weights are always positive, as well as the degree of the product.
  Hence, whenever classical descendant Gromov-Witten invariants are negative
  (e.g. $\langle \tau_1(1) \tau_0(2) \rangle_1^{\text{alg}} = -1$),
  we have an example of classical invariants that do \emph{not} coincide
  with their tropical counterparts as we define them. 
\end{remark}

\section {Splitting curves} \label {splittingcurves}


Now we want to think of curves with a contracted bounded edge as reducible curves. We do that basically by cutting the contracted bounded edge. We have to show that the weight of the curve $C$ is (almost) the product of the two weights of the two curves that arise after cutting. We use the description of weight in terms of determinants from section \ref{multiplicities}.

\begin{notation}
As in section \ref{contractededge} we assume 
$1 \in L$, $2 \in M$ and $3,4 \in N$.
Additionally we require from now on $L = \{1\}$ (i.e.
the marked end $1$ is the only ``free'' end) and $r_3 \geq 1$ (i.e. the marked 
end is restricted by at least one Psi-class).

Let $Z$ be a zero-dimensional cycle of the form $(\ft^*(\lambda) \cdot 
\tau_0(0) \tau_0(1)^m \prod_{k\in N} \tau_{r_k}(2))_d$
with general condition $G_j$ and $P_k$ and $\lambda$, as in \ref{not-z}. 
Let $C\in Z$.
\end{notation}
\begin{construction}\label{cons-cut}
Assume $C$ satisfying $|\det_C(\ev\times \ft)|\neq 0$ has a contracted bounded edge $e$. Cut the bounded edge $e$, thus producing two marked ends. In this way we get two curves $C_1$ and $C_2$ that both have a new marked end in the place of $e$. 
Let $L_i$, $M_i$ and $N_i$ be the subsets of $L$, $M$ and $N$ of marked ends in $C_i$. 
Let $l_i$, $m_i$ and $n_i$ be the sizes of these subsets.
Let $d_i$ be the degree of $C_i$.
Denote by 
\begin{displaymath}
  \ev_{ M_i\cup N_i}:
  \prod_{j\in N_i}\psi_j^{r_j}\cdot[\mathcal{M}_{0,l_i+n_i+m_i+1}^{\lab}(\RR^2,d_i)]
  \rightarrow\RR^{m_i+2n_i}
\end{displaymath} 
the map that evaluates the first coordinate for the points in $M_i$ and both coordinates for the ends in $N_i$ (as in \ref{not-z}). Denote by 
\begin{displaymath}
  \ev_e:\prod_{j\in N_i}\psi_j^{r_j}
  \cdot[\mathcal{M}_{0,l_i+n_i+m_i+1}^{\lab}(\RR^2,d_i)]
  \rightarrow\RR^{2}
\end{displaymath} 
the evaluation at $e$ at both coordinates, and by $(\ev_e)_x$ the evaluation at the first coordinate.
Denote by $\tilde{C_i}$ the curve $C_i$ where we remove the marked end $e$ and straighten the $2$-valent vertex which appears.

Let
\begin{displaymath}
  Z_i:=\prod_{j\in M_i}\ev^{\ast}_{j}(G_j)\cdot \prod_{k\in N_i}\ev^{\ast}_{k}(P_k)\cdot
  \psi_k^{r_k}\cdot[\mathcal{M}_{0,l_i+n_i+m_i+1}^{\lab}(\RR^2,d_i)]
\end{displaymath}
denote the corresponding intersection products. 
\end{construction}

\begin{notation}
We pullback a general point $\lambda\in \calM_{0,4}$ (i.e.\ not the vertex) via the forgetful map $\ft: \mathcal{M}_{0,1+m+n}^{\lab}(\RR^2,d)\rightarrow \calM_{0,4}$. There are $3$ types of such general points, corresponding to the $3$ types of abstract tropical curves with $4$ marked ends. The ends $1$ and $2$ can be together at a vertex, or the ends $1$ and $3$, or the ends $1$ and $4$. We use the following short notation: if $\ft(C)$ is in the ray corresponding to the type where $1$ and $2$ are together at a vertex, we say $\ft(C)=12/34$ (and analogously
in the other cases).
\end{notation}

\begin{lemma}\label{lem-Ccut}
Let $C$ be as in construction \ref{cons-cut} and stick to the notations from there.
If $\ft(C)=12/34$ (then \ $\{1\} = L_1$, so $l_1=1$ and $l_2=0$), then either $d_1=0$ and $L_1\cup M_1\cup N_1=\{1,2\}$, or $d_1,d_2>0$ and there are $3$ cases to distinguish (of which the first and last are symmetric):
\begin{enumerate}
\item\label{item-lem-Ccut-1} $\dim(Z_1)=0$ and $\dim(Z_2)=2$,
\item \label{item-lem-Ccut-2}$\dim(Z_1)=1$ and $\dim(Z_2)=1$, or
\item \label{item-lem-Ccut-3}$\dim(Z_1)=2$ and $\dim(Z_2)=0$.
\end{enumerate}
If $\ft(C)=13/24$ then $d_1,d_2>0$ and the analogous $3$ cases are to distinguish.
\end{lemma}
\begin{proof}
If there were two contracted edges, then all evaluations (i.e.\ $2n+m$ coordinates) would depend only on $1+m+n+3d-3-\sum r_i=m+2n-1$ coordinates, so we get $|\det_C(\ev\times \ft)|= 0$. So we can assume now there is only one contracted bounded edge $e$.
Since $e$ has to count towards the $\mathcal{M}_{0,4}$-coordinate to satisfy $|\det_C(\ev\times \ft)|\neq 0$, $1,j\subset L_1\cup M_1\cup N_1$ and $k,l\subset L_2\cup M_2\cup N_2$ if $\ft(C)=1j/kl$. \\
Let us first consider the case where one of the $d_i$'s is zero. This implies that all edges of the corresponding
curve $C_i$ are contracted. As we cannot have more contracted bounded edges, $C_i$ is a star-shaped curve
containing only a single vertex $V$. But \ref{GeneralValence} states that the ends $2,3,4 \in M \cup N$
all lie at different vertices. Thus the cases $d_1 = 0, \ft(C)=13/24$ and $d_2 = 0$ cannot occur, 
whereas in the remaining case $d_1 = 0, \ft(C)=12/34$ 
the single vertex $V$ must be $3$-valent which is the same as $L_1\cup M_1\cup N_1=\{1,2\}$. \\
Let us now assume $d_1,d_2>0$.
It remains to show that $\dim(Z_1)+\dim(Z_2)=2$ which follows since
\begin{align*}
&\dim(Z_1)+\dim(Z_2)\\=& 3d_1-\sum_{k\in N_1}r_k-m_1-2n_1+3d_2-\sum_{k\in N_2}r_k-m_2-2n_2 \\& = 3d-\sum_{k\in N} r_k-m-2n= 2
\end{align*}
where the last equality follows since $Z$ is zero-dimensional and thus $3d-\sum_{k\in N} r_k-1=m+2n+1$.
\end{proof}

\begin{remark}\label{remark-basis}
In the following, we will choose bases in order to write down an explicit matrix representation for the map $\ev \times \ft$ or $\ev$ locally on a cone.
For a cone $\sigma$ of $\M{n}{d}$ corresponding to a combinatorial type (i.e.\ an abstract graph $\Gamma$ (without length) together with all direction vectors) we pick a root vertex $V$ of $\Gamma$ and choose the coordinates of the point $h(V)\in \RR^2$ to which this vertex is mapped as two coordinates. The remaining coordinates of $\sigma$ are given by the lengths of the bounded edges.
For the spaces $\RR^2$ or $\RR$ that describe our incidence conditions locally, we choose the standard basis vectors. 
It follows from remark 3.2 of \cite{GM053} that the absolute value of the determinant does not depend on any of the choices we make. 
\end{remark}

\begin{lemma}\label{lem-cutmult}
Let $C$ be as in construction \ref{cons-cut} and stick to the notations from there.
If $\ft(C)=12/34$ and $d_1=0$ we want to show $|\det_C(\ev\times \ft)|=(G_2\cdot \tilde{C_2})_{h(e)}\cdot|\det_{\tilde{C_2}}(\ev_{M_2\cup N_2})|$.
For the other three cases from lemma \ref{lem-Ccut} we want to show:
\begin{enumerate}
\item $|\det_C(\ev\times \ft)|=|\det_{C_1}(\ev_{M_1\cup N_1})|\cdot  | \det_{C_2}(\ev_{M_2\cup N_2}\times \ev_e)|$,
\item $|\det_C(\ev\times \ft)|=(\tilde{C_1}\cdot \tilde{C_2})_{h(e)}\cdot|\det_{\tilde{C_1}}(\ev_{M_1\cup N_1})|\cdot |\det_{\tilde{C_2}}(\ev_{M_2\cup N_2})|$, or
\item $|\det_C(\ev\times \ft)|=|\det_{C_1}(\ev_{M_1\cup N_1}\times \ev_e)|\cdot |\det_{C_2}(\ev_{M_2\cup N_2})|$.
\end{enumerate}
\end{lemma}

\begin{proof}
For all cases, note first that the matrix of $|\det_C(\ev\times \ft)|$ has a column with only zeros except one $1$. This is the column corresponding to $e$. Since $e$ is contracted, it is not needed for any evaluation. But it is needed for the $\mathcal{M}_{0,4}$-coordinate, so it has zeros except a $1$ in the $\ft$-row.
We can delete this row and column without changing the absolute value of the determinant. Call the matrix with the deleted row and column $A$. Then $|\det(A)|=|\det_C(\ev\times \ft)|$.

Now let $\ft(C)=12/34$ and $d_1=0$, it follows $L_1\cup M_1\cup N_1=\{1,2\}$.
We want to show that the boundary vertex $V$ of $e$ in $C_2$ is $3$-valent, too. Assume it is not, then there has to be a marked end with a Psi-condition adjacent to $V$. But this marked end is in $N$ and thus required to meet a point. This is a contradiction, since the point is not on the line that $2$ is required to meet (cf. \ref{GeneralConditions} (d)). So let $ e_1 $ and $ e_2 $ be the two other
  edges adjacent to $V$ and assume first that both of them are bounded. Denote their
  common direction vector (up to sign) by $ v=(v_1,v_2) $ and their lengths by $ l(e_1), l(e_2) $.
  Assume that the root vertex is on the $ e_1 $-side of $ e $. Then the
  entries of the matrix $A$ corresponding to $ l(e_1) $ and $ l(e_2) $ are
    \[ \begin {array}{l|cc}
         \mbox {$ \downarrow $ evaluation at\dots} & l(e_1) & l(e_2) \\ \hline
         \mbox {$ 2 $ (1 row)} & v_1 & 0 \\
         \mbox {points reached via $ e_1 $ from $ 2 $ \mbox{(1 or 2 rows)}}
           & 0 & 0 \\
         \mbox {points reached via $ e_2 $ from $ 2 $ \mbox{(1 or 2 rows)}}
           & v & v
       \end {array} \]
  We see that after subtracting the $ l(e_2) $-column from the $ l(e_1) $-column we
  again get one column with only one non-zero entry $ v_1 $. So for the
  determinant we get $ v_1$ as a factor, dropping the
  corresponding row and column (which means removing $ e$
  and straightening the $2$-valent vertex), so we get $|\det(A)|=v_1\cdot|\det_{\tilde{C_2}}(\ev_{M_2\cup N_2}) | =(\tilde{C_2} \cdot G_2)_{h(e)}\cdot |\det_{\tilde{C_2}}(\ev_{M_2\cup N_2})|$. Essentially the same
  argument holds if one of the adjacent edges --- say $ e_2 $ -- is unbounded:
  in this case there is only an $ l(e_1) $-column which has zeroes everywhere
  except in the one $ 2 $-row where the entry is $ v_1 $.

 Next, let $\dim(Z_1)=0$ and $\dim(Z_2)=2$. 
Denote by $a_i$ the dimension of $$\prod_{k\in N_i}\psi_k^{r_k}\cdot \mathcal{M}^{\lab}_{0+l_i+m_i+n_i+1}(\RR^2,d_i),$$that is, 
$$a_i= 3d_i+l_i+m_i+n_i+1-\sum_{k\in N_i}r_k-1.$$ Since $\dim(Z_1)=0$ we have $m_1+2n_1=a_1$ and since $\dim(Z_2)=2$ we have $m_2+2n_2=a_2-2$.
Let the boundary vertex $V$ of $e$ in $C_1$ be the root vertex for $C$. Choose the following order of coordinates: start with the root vertex, then bounded edges in $C_1$, next bounded edges in $C_2$. 
Start with the marked ends in $C_1$ and then add the marked ends in $C_2$. Then the matrix $A$ is in block form: because the points on $C_1$ need only the root vertex and the bounded edges of $C_1$, they need the first $a_{1}=m_1+2n_1$ coordinates, and have $0$ after that. So there is a $0$ block on the top right, and the top left is just the matrix of $\ev_{M_1\cup N_1}$ at $C_1$. So $|\det(A)|=|\det_{C_1}(\ev_{M_1\cup N_1})|\cdot |\det(B)|$ where $B$ denotes the lower right box. 
Consider the matrix of $\ev_{M_2\cup N_2}\times \ev_e$ at $C_2$, and let the root vertex be the boundary vertex of $e$ in $C_2$. Then this matrix has two more rows and columns than $B$, namely the root vertex columns and the rows corresponding to $\ev_e$. But since these two rows start with a $2\times 2$ unit matrix block and have zeros after that, we can see that $|\det(B)|=|\det_{C_2}(\ev_{M_2\cup N_2}\times \ev_e)|$.

The third case is symmetric.
Finally, assume $\dim(Z_1)=1$ and $\dim(Z_2)=1$, i.e.\ $m_1+2n_1= a_1-1$ and $m_2+2n_2= a_2-1$.
First we want to show that the two vertices of $e$ are $3$-valent. Assume the vertex in $C_1$, $V$, is not $3$-valent, then there must be a marked end $i$ with a Psi-class adjacent to $V$. But this end is in $N$ then, so it is required to meet a point $P_i\in\RR^2$. Since $\dim(Z_1)=1$ we can move $C_1$ locally in a $1$-dimensional family such that all incidence conditions are still defined. Let $C_1'$ be an element of this family. Since $C_1'$ has to meet $P_i$ as well, we can glue $C_1'$ to $C_2$ thus producing a curve $C$ in $Z$. This is a contradiction since the dimension of $Z$ is $0$.

Since the argument is symmetric it follows that both vertices of $e$ are $3$-valent.
Denote the two edges adjacent to $e$ in $C_1$ by $e_1$ and $e_2$ and the two edges in $C_2$ adjacent to $e$ by $e_3$ and $e_4$. Assume first that all of those edges are bounded.
Let the boundary vertex $V$ of $e$ in $C_1$ be the root vertex for $C$. Then the matrix $A$ reads:
    \[ \begin {array}{c|l|ccc@{\;\,}c@{\;\,}c@{\;\,}cc}
         & & & \mbox {lengths in $ C_1 $} & & & & & 
           \mbox {lengths in $ C_2 $} \\
         & & \mbox {root} &
           \mbox {($ a_1-4 $ cols)} & l(e_1) & l(e_2) & l(e_3) & l(e_4) &
           \mbox {($ a_2-4 $ cols)} \\ \hline \hline
       (2n_1+m_1 &
         \mbox {ends behind $ e_1 $} &
        I_2 & * & v & 0 & 0 & 0 & 0 \\
      \mbox {rows}) & \mbox {ends behind $ e_2 $} &
         I_2 & * & 0 & -v & 0 & 0 & 0 \\ \hline
       (2n_2+m_2 &
         \mbox {ends behind $ e_3 $} &
         I_2 & 0 & 0 & 0 & w & 0 & * \\
       \mbox {rows}) & \mbox {ends behind $ e_4 $} &
         I_2 & 0 & 0 & 0 & 0 & -w & *
       \end {array} \]
  where $ I_2 $ is the $ 2 \times 2 $ unit matrix, and $*$ denotes arbitrary entries.
  Now add $v$ times the root columns to the $ l(e_2) $-column, subtract the $ l(e_1)
  $-column from the $ l(e_2) $-column and the $ l(e_4) $-column from the $ l(e_3)
  $-column to obtain the following matrix with the same determinant:
    \[ \begin {array}{c|l|ccc|c@{\;\,}c@{\;\,}cc}
         & & & \mbox {lengths in $ C_1 $} & & & & & 
           \mbox {lengths in $ C_2 $} \\
         & & \mbox {root} &
           \mbox {($ a_1-4 $ cols)} & l(e_1) & l(e_2) & l(e_3) & l(e_4) &
           \mbox {($ a_2-4 $ cols)} \\ \hline \hline
       (2n_1+m_1 &
         \mbox {ends behind $ E_1 $} &
         I_2 & * & v & 0 & 0 & 0 & 0 \\
       \mbox {rows}) & \mbox {ends behind $ E_2 $} &
         I_2 & * & 0 & 0 & 0 & 0 & 0 \\ \hline
       (2n_2+m_2 &
         \mbox {ends behind $ E_3 $} &
         I_2 & 0 & 0 & v & w & 0 & * \\
       \mbox {rows}) & \mbox {ends behind $ E_4 $} &
         I_2 & 0 & 0 & v & w & -w & *
       \end {array} \]
  Note that this matrix has a block form with a zero block at the top right.
  Denote the top left block (of size $ 2n_1+m_1=2+a_1-4+1 $) by $ A_1 $ and the bottom
  right (of size $ 2n_2+m_2=3+a_2-4$) by $ A_2 $, then $|\det(A)|= | \det A_1 \cdot \det A_2 | $.

  The matrix $ A_1 $ is precisely the matrix for the evaluation map $\ev_{M_1\cup N_1}$ of $ \tilde{C_1} $ (which arises from $C_1$ after forgetting the marked end corresponding to $e$) if we choose the other vertex
  of $ e_2 $ as the root vertex. Hence $ |\det A_1| = |\det_{\tilde{C_1}}(\ev_{M_1\cup N_1})|$.
  In the same way the matrix for the evaluation map $\ev_{M_2\cup N_2}$ of $ \tilde{C_2} $, if we again
  forget the marked end corresponding to $e$ and now choose the other vertex
   of $ e_3 $ as the root vertex, is the matrix $ A_2' $ obtained from $
  A_2 $ by replacing $v$ and $w$ in the first two columns by the first and
  second unit vector, respectively. But $ A_{2} $ is simply obtained from $
  A_2' $ by right multiplication with the matrix
    \[ \left( \begin {array}{ccc}
         v & w & 0 \\
         0 & 0 & I_{2n_2+2}
       \end {array} \right) \]
  which has determinant $ \det (v,w) $. So we conclude that
    \begin{displaymath} |\det A_2| = |\det (v,w)| \cdot |\det A_2'|
         = (C_1 \cdot C_2)_{h(e)} \cdot|\mbox{$\det _{\tilde{C_2}}(\ev_{M_2\cup N_2})$}|. \end{displaymath}
\end{proof}

\begin{remark}\label{rem-converse}
 The following ``converse'' of lemma \ref{lem-Ccut} and lemma \ref{lem-cutmult} is also
  true: 
For each choice of $\tilde{C_2}$ satisfying all conditions but $2$ and each choice of an intersection point of $\tilde{C_2}$ with $G_2$ we can add a contracted bounded edge and the two marked ends $1,2$ on the other side to built exactly one possible $C$. The curve $C$ then contributes $(G_2\cdot \tilde{C_2})_{h(e)}|\det_{\tilde{C_2}}(\ev_{M_2\cup N_2})|$ to the count. By B\'ezout's theorem (\cite{RST03}), each choice of $\tilde{C_2}$ contributes $d_2\cdot |\det_{\tilde{C_2}}(\ev_{M_2\cup N_2})|=d\cdot |\det_{\tilde{C_2}}(\ev_{M_2\cup N_2})|$.

For each choice of $C_1$ satisfying the conditions in $L_1\cup M_1\cup N_1$ and each choice of $C_2$ satisfying the conditions in $L_2\cup M_2\cup N_2$ plus in addition the condition $h(e)=p$ we get exactly one possible $C$ by gluing the two curves along $e$. This curve $C$ contributes to the count with weight $|\det_{C_1}(\ev_{M_1\cup N_1})|\cdot |\det_{C_2}(\ev_{M_2\cup N_2}\times \ev_e)|$ (and the other way round).

For each choice of $\tilde{C_1}$ and $\tilde{C_2}$ satisfying the conditions in $L_1\cup M_1\cup N_1$ and $L_2\cup M_2\cup N_2$ and for each choice of points $ P \in \tilde{C_1} $ and $ Q \in \tilde{C_2} $ that map to
  the same image point in $\RR^2 $ we can glue $P$ and $Q$ along a contracted bounded edge and thus built exactly one possible $C$. The curve $C$ contributes to the count with weight $(C_1\cdot C_2)_{h(e)}\cdot|\det_{\tilde{C_1}}(\ev_{M_1\cup N_1})|\cdot |\det_{\tilde{C_2}}(\ev_{M_2\cup N_2})|$.
By B\'ezout's theorem, each choice of $\tilde{C_1}$ and $\tilde{C_2}$ thus contributes $d_1\cdot d_2\cdot |\det_{\tilde{C_1}}(\ev_{M_1\cup N_1})|\cdot| \det_{\tilde{C_2}}(\ev_{M_2\cup N_2})|$.

\end {remark}

\section {String and divisor equation} \label {stringdivisor}


In this section we prove two lemmas which deal
with the case of an extra end in a top-dimensional intersection product that
is restricted either by no condition at all (string equation) or by only 
a line condition (divisor equation). 

\begin{lemma}[String equation]
For tropical descendant Gromov-Witten invariants the following equality holds:
\begin{align*}
  & \langle \tau_0(0) \cdot\tau_0(0)^l\cdot\tau_0(1)^m\cdot 
    \prod_{k\in N}\tau_{r_k}(2)\rangle_d \\
  = & \sum_{\twolines{k\in N}{r_k > 0}}\langle \tau_0(0)^l\cdot\tau_0(1)^m \cdot 
    \tau_{r_k-1}(2)\cdot \prod_{k\neq k'\in N}\tau_{r_{k'}}(2)\rangle_d
\end{align*}
\end{lemma}

\begin{proof}
  Choose incidence conditions $G_j, P_k$ 
  such that they are general for all the derived intersection products
  \begin{displaymath}
    Z := (\tau_0(0) \cdot\tau_0(0)^l\cdot\tau_0(1)^m\cdot 
  \prod_{k\in N}\tau_{r_k}(2))_d,
  \end{displaymath} 
  \begin{displaymath}
    Z_k := (\tau_0(0)^l\cdot\tau_0(1)^m \cdot 
    \tau_{r_k-1}(2)\cdot \prod_{k\neq k'\in N}\tau_{r_{k'}}(2))_d
  \end{displaymath} 
  (note that
  $Z$ lives in $\M{1+l+m+n}{d}$, whereas the $Z_k$ lives in $\M{l+m+n}{d}$).
  Then \ref{GeneralConditions} (a) tells us that the products
  just consist of the set of curves fulfilling the incidence conditions and
  having required valences, with the additional data of a weight
  for each curve. \\
  Let $C'$ be a curve in $Z_k$. Then we obtain a curve $C \in Z$ by
  attaching the additional end, say $x$, to the vertex $V_k$ at which the end $k$ lies.
  Let us check that the weight of $C'$ in $Z_k$ and $C$ in $Z$ coincide.
  As our conditions are general, $C'$ lies in a facet $\sigma'$ of 
  $\psi_k^{r_k - 1}\cdot \prod_{k\neq k'\in N}\psi_{k'}^{r_{k'}}
  \cdot \M{l+m+n}{d}$ and $C$ lies in a facet 
  $\sigma$ of $\prod_{k\in N}\psi_k^{r_k} \cdot \M{1+l+m+n}{d})$. Moreover,
  the map 
  \begin{displaymath}
    \ft_x : \M{1+l+m+n}{d} \rightarrow \M{l+m+n}{d}
  \end{displaymath}
  forgetting
  the additional end $x$ maps $\sigma$ $\ZZ$-isomorphically to $\sigma'$ (the
  inverse is given by adding $1$ to $V_k$ as above). The evaluation maps
  $\ev_k$ on $\sigma$ are just obtained as pull-backs $\ft_{x*}(\ev'_k)$,
  where $\ev'_k$ denotes the corresponding evaluation map on $\sigma'$.
  Hence, the weights of $C$ and $C'$ coincide. \\
  It remains to check that each $C \in Z$ is obtained in the above way
  from $C' \in Z_k$ for unique $k \in N$. Uniqueness is clear, as by
  \ref{GeneralValence} (b') all ends $k \in N$ lie at pairwise different
  vertices and hence $x$ can not be adjacent to more than one end $k \in N$.
  On the other hand, to show that it is adjacent to a $k \in N$ with $r_k > 0$,
  it suffices to show that $x$ cannot be adjacent to a $3$-valent vertex.
  If it were, at least one of the other two adjacent edges, say $E$ would be bounded
  (otherwise the abstract graph were not connected). But then, we could
  change the length of $E$ (and accordingly the length of the other edge
  if necessary) without changing the coordinates of the marked ends, which
  contradicts the fact that the set of curves fulfilling our given conditions
  is finite by \ref{GeneralConditions} (a). \\
  \end{proof}

\begin{lemma}[Divisor equation] \label{DivisorEquation}
For tropical descendant Gromov-Witten invariants the following equality holds:

\begin{align*}\langle \tau_0(1) \cdot\tau_0(0)^l\cdot\tau_0(1)^m\cdot \prod_{k\in N}\tau_{r_k}(2)\rangle_d
=d\cdot\langle \tau_0(0)^l\cdot\tau_0(1)^m\cdot \prod_{k\in N}\tau_{r_{k}}(2)\rangle_d\end{align*}

\end{lemma}

\begin{proof}
First we choose general incidence conditions.
Because of \ref{InvariantsAndIndependence} we can assume that the line conditions are all vertical lines, i.e.\ of type $\max\{x,0\}$.
Then for all curves $C$ in $( \tau_0(1) \cdot\tau_0(0)^l\cdot\tau_0(1)^m\cdot 
\prod_{k\in N}\tau_{r_k}(2))_d$ we know that their 
weight is equal to $|\det_{C}(\ev)|$, where $\ev$ denotes the product of all evaluation maps (evaluation of the $x$-coordinate for all lines, both coordinates for all points) (theorem \ref{MultEqualDet}).
Assume $x$ is the additional marked end with line condition $G$ (but 
without Psi-condition).
$x$ has to be adjacent to a $3$-valent vertex (see \ref{GeneralValence}).
Exactly as in lemma \ref{lem-cutmult} we can see that $|\det_{C}(\ev)|
=(G\cdot \tilde{C})_{h(x)}\cdot|\det_{\tilde{C}}(\ev_x)|$ where 
$\tilde{C}$ is the curve we get when forgetting $x$ (i.e.\ removing it from $C$ and 
straightening the $2$-valent vertex) and $\ev_x$ is the product of all other evaluations.
Thus any curve in $( \tau_0(1) \cdot\tau_0(0)^l\cdot\tau_0(1)^m
\cdot \prod_{k\in N}\tau_{r_k}(2))_d$ gives us a curve in 
$(\tau_0(0)^l\cdot\tau_0(1)^m\cdot \prod_{k\in N}\tau_{r_{k}}(2)
)_d$ by removing the marked end $x$. Conversely, given a 
curve $\tilde{C}$ in $( \tau_0(0)^l\cdot\tau_0(1)^m\cdot 
\prod_{k\in N}\tau_{r_{k}}(2))_d$ we can pick a point 
$p\in (G\cdot \tilde{C})$ and attach a marked end to get a curve 
$C\in ( \tau_0(1) \cdot\tau_0(0)^l\cdot\tau_0(1)^m\cdot \prod_{k\in N}
\tau_{r_k}(2))_d$. Since $(G\cdot \tilde{C})=\sum_p 
(G\cdot \tilde{C})_p=d$ by tropical B\'ezout's theorem (\cite{RST03}), the statement follows.
\end{proof}

\section {Recursion} \label {recursion}


Now we sum up the results of the preceding sections to a certain WDVV equation.
We also show in this section that this WDVV equation together with the string and the divisor equation
are sufficient to show that the tropical invariants coincide
with the classical ones.
To distinguish our tropical invariants that we denote by $\langle\tau_0(0)^l \tau_0(1)^m \prod_{k\in N} \tau_{r_k}(2)\rangle_d$ from the classical ones, we use the notation $\langle\tau_0(0)^l \tau_0(1)^m \prod_{k\in N} \tau_{r_k}(2)\rangle_d^{\text{alg}}$ for the classical invariants.

\begin{theorem}\label{thm-WDVV}
The tropical descendant invariants as defined in \ref{InvariantsAndIndependence} satisfy the following WDVV equation if $r_3 > 0$:
\begin{align*}
& \langle  \tau_0(1)^{m} \prod_{k\in N} \tau_{r_k}(2) \rangle_{d}+
\\& \sum D\cdot\langle\tau_0(0) \tau_0(1) \tau_0(1)^{m_1} \prod_{k\in N_1} \tau_{r_k}(2)\tau_0(e)\rangle_{d_1} \cdot \langle \tau_{r_3}(2)\tau_{r_4}(2) \tau_0(1)^{m_2} \prod_{k\in N_2} \tau_{r_k}(2)\tau_0(f)\rangle_{d_2}
\\=&\sum D\cdot\langle\tau_0(0) \tau_{r_3} (2) \tau_0(1)^{m_1} \prod_{k\in N_1} \tau_{r_k}(2) \tau_0(e)\rangle_{d_1} \cdot \langle\tau_0(1) \tau_{r_4}(2) \tau_0(1)^{m_2} \prod_{k\in N_2} \tau_{r_k}(2)\tau_0(f)\rangle_{d_2}
\end{align*}
where $$D= \frac{(d_1!)^3\cdot (d_2!)^3}{d!^3}$$ and the sums range over all 
\begin{align*}&e+f=2, \;e,f \geq 0, \\&M_1\dunion M_2=M\setminus\{2\}, \\&N_1\dunion N_2= N\setminus\{3,4\}\mbox{ and }\\&E_1\dunion E_2=E, E_1,E_2\neq \emptyset.\end{align*} Here, $E$ denotes the set of non-contracted ends, and $E_1$ is subset of non-contracted ends such that each of the standard directions $-e_1,-e_2,e_1+e_2$ appears $d_1$ times. 

The equation can be rewritten as
\begin{align*}
& \langle  \tau_0(1)^{m} \prod_{k\in N} \tau_{r_k}(2) \rangle_{d}+
\\& \sum \langle\tau_0(0) \tau_0(1) \tau_0(1)^{m_1} \prod_{k\in N_1} \tau_{r_k}(2)\tau_0(e)\rangle_{d_1} \cdot \langle \tau_{r_3}(2)\tau_{r_4}(2) \tau_0(1)^{m_2} \prod_{k\in N_2} \tau_{r_k}(2)\tau_0(f)\rangle_{d_2}
\\=&\sum \langle\tau_0(0) \tau_{r_3} (2) \tau_0(1)^{m_1} \prod_{k\in N_1} \tau_{r_k}(2) \tau_0(e)\rangle_{d_1} \cdot \langle\tau_0(1) \tau_{r_4}(2) \tau_0(1)^{m_2} \prod_{k\in N_2} \tau_{r_k}(2)\tau_0(f)\rangle_{d_2}
\end{align*}
where now the sums range over all 
\begin{align*}&e+f=2, \\&M_1\dunion M_2=M\setminus\{2\}, \\&N_1\dunion N_2= N\setminus\{3,4\}\mbox{ and }\\&d_1+d_2=d, d_1,d_2>0.\end{align*}

\end{theorem}
\begin{proof}
It follows from rational equivalence that
\begin{align*} &\langle\tau_0(0) \tau_0(1)^{m} \prod_{k\in N} \tau_{r_k}(2) \ft^*(\lambda)\rangle_{d} \\&=
\frac{1}{(d!)^3}\deg\left( \prod_{j\in M}\ev_j^*(G_j)\prod_{k\in N} \ev_k^*(P_k)\psi_k^{r_k}\cdot \ft^*(\lambda)\cdot \M{1+m+n}{d}\right)\end{align*} does not depend on the choice of $\lambda\in \calM_{0,4}$. Thus we can pick a very large $\lambda_1$ on the ray $12/34$ of $\calM_{0,4}$ and a very large $\lambda_2$ on the ray $13/24$ and set the degree equal for those two values.
Denote by 
$$Z:= \deg\left( \prod_{j\in M}\ev_j^*(G_j)\prod_{k\in N} \ev_k^*(P_k)\psi_k^{r_k}\cdot \ft^*(\lambda_1)\cdot \M{1+m+n}{d}\right).$$
We show that the left hand side of the above sum equals $\frac{1}{d!^3}$ times the degree of $Z$. Analogously one can show that the right hand side equals $\frac{1}{d!^3}$ times the degree of the analogous intersection product with $\lambda_2$, which finishes the proof.
By theorem \ref{MultEqualDet} we know that 
 $$
    Z = \sum_{C \in S} |{\det}_C(\ft \times \ev)| \cdot C,
  $$  
where $S$ is the set of curves in $\M{1+m+n}{d} $ satisfying the point and line conditions and mapping to $\lambda_1$ under $\ft$.
Let $F=(\tau_0(0) \tau_0(1)^{m} \prod_{k\in N} \tau_{r_k}(2))_{d}$, then $F$ is a one-dimensional cycle.
Let $\sigma$ be a cone of $F$ corresponding to curves without a contracted bounded edge. Then lemma \ref{ContractedEdge} says that the image of $\sigma$ under $\ft$ is bounded. Since we picked $\lambda_1$ to be very large, we therefore know that $\sigma$ cannot contribute to the degree of $Z$. Hence all $C\in S$ contain a contracted bounded edge.
Pick a curve $C\in S$, then we know by \ref{lem-Ccut} that we can cut the contracted edge thus producing two curves $C_1$ and $C_2$ with an extra marked end $e$.

If the degree of $C_1$, $d_1$, equals $0$ then we know by \ref{lem-cutmult}
that $$|{\det}_C(\ft \times \ev)|= (G_2\cdot \tilde{C_2})_{h(e)}\cdot|{\det}_{\tilde{C_2}}(\ev_{M_2\cup N_2})|,$$ where $G_2$ denotes the line condition for the marked end $2$ and $\tilde{C_2}$ denotes the curve that we get from $C_2$ by forgetting the additional marked end $e$.
By \ref{rem-converse} we know that each choice of $\tilde{C_2}$ satisfying all conditions in $L_2\cup M_2\cup N_2=L\cup M\cup N\setminus\{1,2\}$ contributes $d\cdot |{\det}_{\tilde{C_2}}(\ev_{M_2\cup N_2})|$ possible curves $C$ (counted with weight).
Thus the contribution to $Z$ from curves $C$ such that $d_1=0$ equals
$$ d\cdot \prod_{j\in M\setminus\{2\}}\ev_j^*(G_j)\prod_{k\in N}\ev_k^*(P_k)\psi_k^{r_k}\cdot \M{m-1+n}{d} $$
which by the divisor equation (\ref{DivisorEquation}) equals
$$  \prod_{j\in M}\ev_j^*(G_j)\prod_{k\in N}\ev_k^*(P_k)\psi_k^{r_k}\cdot \M{m+n}{d} .$$
Multiplying by the factor $\frac{1}{d!^3}$, we can see that those curves contribute
$$\langle  \tau_0(1)^{m} \prod_{k\in N} \tau_{r_k}(2) \rangle_{d}$$
to $\frac{1}{d!^3} \deg Z$.

Now assume that $d_1>0$ and denote as in \ref{cons-cut} 
$$Z_i:=\prod_{j\in M_i}\ev^{\ast}_{j}(G_j)\cdot \prod_{j\in N_i}\ev^{\ast}_{j}(p_j)\cdot \prod_{j\in N_i}\psi_j^{r_j}\cdot[\mathcal{M}_{0,l_i+n_i+m_i+1}^{\lab}(\RR^2,d_i)].$$

Then we know by \ref{lem-Ccut} that one of the following three cases hold:
\begin{enumerate}
\item $\dim(Z_1)=0$ and $\dim(Z_2)=2$ or
\item $\dim(Z_1)=1$ and $\dim(Z_2)=1$ or
\item $\dim(Z_1)=2$ and $\dim(Z_2)=0$.
\end{enumerate}
We know by \ref{lem-cutmult} that in the first case,
$$|{\det}_C(\ev\times \ft)|=|{\det}_{C_1}(\ev_{M_1\cup N_1})|\cdot  | {\det}_{C_2}(\ev_{M_2\cup N_2}\times \ev_e)|,$$
where $\ev_e$ now denotes the evaluation on both coordinates of the new marked end $e$.
By \ref{rem-converse} we know that for each choice of $C_1$ and $C_2$ satisfying the conditions we get exactly one possible $C$.
But by \ref{MultEqualDet} we know that 
 $$
    Z_1 = \sum_{C_1} |{{\det}}_{C_1}( \ev_{M_1\cup N_1})| \cdot C_1,
  $$  
and analogously
 $$
   \ev_e^*(P)\cdot Z_2 = \sum_{C_2} |{{\det}}_{C_2}( \ev_{L_2\cup M_2\cup N_2}\times \ev_e)| \cdot C_2,
  $$

Thus we get a contribution of $\deg(Z_1)\cdot \deg( \ev_e^*(P)\cdot Z_2)$ to $\deg(Z)$, respectively
$$ \frac{(d_1!)^3\cdot (d_2!)^3}{d!^3} \langle\tau_0(0) \tau_0(1) \tau_0(1)^{m_1} \prod_{k\in N_1} \tau_{r_k}(2)\rangle_{d_1} \cdot \langle \tau_{r_3}(2)\tau_{r_4}(2) \tau_0(1)^{m_2} \prod_{k\in N_2} \tau_{r_k}(2)\tau_0(2)\rangle_{d_2}$$ to $\frac{1}{d!^3} \deg Z$.

Analogously, we get a contribution of 
$$ \frac{(d_1!)^3\cdot (d_2!)^3}{d!^3} \langle\tau_0(0) \tau_0(1) \tau_0(1)^{m_1} \prod_{k\in N_1} \tau_{r_k}(2)\tau_0(2)\rangle_{d_1} \cdot \langle \tau_{r_3}(2)\tau_{r_4}(2) \tau_0(1)^{m_2} \prod_{k\in N_2} \tau_{r_k}(2)\rangle_{d_2}$$
in the last case.

In the second case, we know by \ref{lem-cutmult} that 
$$|{\det}_C(\ev\times \ft)|=(\tilde{C_1}\cdot \tilde{C_2})_{h(e)}\cdot|{\det}_{\tilde{C_1}}(\ev_{M_1\cup N_1})|\cdot |{\det}_{\tilde{C_2}}(\ev_{M_2\cup N_2})|$$

and by \ref{rem-converse} we know that each choice of $\tilde{C_1}$ and $\tilde{C_2}$ satisfying the conditions gives us
$$d_1\cdot d_2\cdot |{\det}_{\tilde{C_1}}(\ev_{M_1\cup N_1})|\cdot| {\det}_{\tilde{C_2}}(\ev_{M_2\cup N_2})|.$$

Since

$$
(\ft_e)_*( Z_i) = \sum_{\tilde{C_i}} |{{\det}}_{\tilde{C_i}}( \ev_{ M_i\cup N_i})| \cdot \tilde{C_i}$$
(where $\ft_e$ denotes the map which forgets the marked point $e$) and since
$$d_1\cdot \langle\tau_0(0) \tau_0(1) \tau_0(1)^{m_1} \prod_{k\in N_1} \tau_{r_k}(2)\rangle_{d_1}
=  \langle\tau_0(0) \tau_0(1) \tau_0(1)^{m_1} \prod_{k\in N_1} \tau_{r_k}(2)\tau_0(1)\rangle_{d_1}$$
and
$$ d_2\cdot  \langle \tau_{r_3}(2)\tau_{r_4}(2) \tau_0(1)^{m_2} \prod_{k\in N_2} \tau_{r_k}(2)\rangle_{d_2}
=\langle \tau_{r_3}(2)\tau_{r_4}(2) \tau_0(1)^{m_2} \prod_{k\in N_2} \tau_{r_k}(2)\tau_0(1)\rangle_{d_2}
$$
by the divisor equation we get a contribution of
$$ \frac{(d_1!)^3\cdot (d_2!)^3}{d!^3} \langle\tau_0(0) \tau_0(1) \tau_0(1)^{m_1} \prod_{k\in N_1} \tau_{r_k}(2)\tau_0(1)\rangle_{d_1} \cdot \langle \tau_{r_3}(2)\tau_{r_4}(2) \tau_0(1)^{m_2} \prod_{k\in N_2} \tau_{r_k}(2)\tau_0(1)\rangle_{d_2}.$$

Finally, there are $\binom{d}{d_1}^3$ choices of the sets $E_1$ and $E_2$ if we fix $d_{1}$ and $d_2$.

\end{proof}

\begin{lemma}
  Choose strictly positive integers $r,d$ such that $1 + 3d - 3 + 2 = 2 + r$.
  Then the classical one-marked-point invariant $ \langle \tau_r(2) \rangle _d^{\text{alg}}$ equals
  $$
    \langle \tau_r(2) \rangle _d^{\text{alg}} = \frac{1}{(d!)^3}.
  $$
\end{lemma}

\begin{proof}
  We use two (classical) WDVV equations (\cite{FP}
  or, more detailed but unpublished, \cite{K}) with four marked points. 
If we compute $\langle \tau_0(0)\tau_0(1)^2\tau_r(2)\ft^*(\lambda)\rangle_d ^{\text{alg}}$ for the two special points $\lambda=12/34$ and $\lambda=13/24$ on $M_{0,4}$, then we get 
$$
    \langle \tau_0(1)^2 \tau_{r}(2) \rangle _d^{\text{alg}} 
      = \langle \tau_0(0) \tau_{r}(2) \tau_0(2) \rangle _d^{\text{alg}}
  $$
as illustrated by the following picture:
  \begin{center}
    \input{WDVV1.pstex_t}
  \end{center}  
 For $\langle \tau_0(0)\tau_0(1)\tau_0(2)\tau_{r-1}(2)\ft^*(\lambda)\rangle_d ^{\text{alg}}$
we get
 \begin{eqnarray*}
    \langle \tau_0(1) \tau_{r-1}(2) \tau_0(2) \rangle _d^{\text{alg}} 
      = \langle \tau_0(0)^2 \tau_{r-1}(2) \rangle _{d-1}^{\text{alg}}
  \end{eqnarray*}
as illustrated by
  \begin{center}
    \input{WDVV2.pstex_t}
  \end{center}
  Now, applying string and divisor equation where possible and plugging in the 
  left hand side of the first equation in the right hand side of the second 
  equation produces
  $$
    d^3 \cdot \langle \tau_r(2) \rangle _d^{\text{alg}} 
      = \langle \tau_{r-3}(2) \rangle _{d-1}^{\text{alg}}.
  $$
  Together with the initial invariant 
  $\langle \tau_0(0)^2 \tau_{0}(2) \rangle _0^{\text{alg}}
  = 1$, this proves the lemma.
\end{proof}

\begin{lemma}
  Choose strictly positive integers $r,d$ such that $1 + 3d - 3 + 2 = 2 + r$.
  Then the tropical one-marked-end invariant $\langle \tau_r(2) \rangle _d^{\text{trop}}$ equals
  $$
    \langle \tau_r(2) \rangle _d^{\text{trop}} = \frac{1}{(d!)^3}.
  $$
\end{lemma}

\begin{proof}
  Choosing the single marked end $e$ as root vertex,
  we get $\M{1}{d} = \calM_{0,1+3d} \times \RR^2$ and the two projections
  are $\ft'$ and $\ev_e$. Recall that Psi-classes
  for parameterized curves are just pull-backs of Psi-classes along $\ft$.
  Using \cite[9.6]{AR07}, we get
  \begin{eqnarray*}
     \langle \tau_r(2) \rangle _d^{\text{trop}} & = &
      \frac{1}{(d!)^3} \deg\left(
      (\ft'(\psi_e))^r \cdot \ev_e(P) \cdot (\calM_{0,1+3d} \times \RR^2)
      \right) \\
    & = &
      \frac{1}{(d!)^3} \deg\left(
      (\psi_e^{r} \cdot \calM_{0,1+3d}) \times (P \cdot \RR^2)
      \right) \\
    & = & \frac{1}{(d!)^3},        
  \end{eqnarray*}
  where in the last step we use $\deg(\psi_e^{r} \cdot \calM_{0,1+3d}) = 1$
  (cf. \cite[4.2]{KM07}).
\end{proof}

\begin{theorem}
  Let $d,l,m,n$ and $r_k, k \in N$ be positive integers with $d>0$ 
  such that
  $$
    l + m + n + 3d - 3 + 2 = m + 2n + \sum_{k \in N} r_k.
  $$
  Then the corresponding tropical and classical descendant invariants coincide,
  i.e.
  $$
     \langle\tau_0(0)^l\tau_0(1)^m \prod_{k\in N}\tau_{r_k}(2) \rangle_d^{\text{trop}} =
      \langle\tau_0(0)^l\tau_0(1)^m \prod_{k\in N}\tau_{r_k}(2) \rangle_d^{\text{alg}} .
  $$
\end{theorem}

\begin{proof}
  The tropical string, divisor and WDVV equations proved in the preceding sections 
  are also fulfilled by the corresponding classical invariants. Hence, we can use
  these equations to reduce our invariants to such ones for which we know or can prove
  that they coincide. \\
  1. case: $r_k = 0$ for all $k \in N$ (i.e. no Psi-classes) \\
  After applying string and divisor equation, we can assume
  that $l = 0 = m$. Using \ref{MultEqualDet}, we see that
  the numbers $I^{\text{trop}}(d;0,0,n;0)$
  are equal to the numbers $N_d$ defined in \cite[3.4, 3.9]{GM053}.
  It is well-known that these numbers coincide
  with the classical ones (see \cite{Mi03} and \cite[5.6]{GM053}). \\
  2. case: exists $k \in N$ with $r_k > 0$ (i.e. at least one Psi-class) \\
  subcase I: $n = 1$ \\
  After applying string and divisor equation, we can assume
  that $l = 0$ and $m = 0$. The two last preceding lemmas show
  that in this case the classical and tropical invariants 
  coincide. \\  
  subcase II: $n \geq 2$ \\
  After applying string and divisor equation, we can assume
  that $l = 0$ and $m = 1$. In particular, if $m=0$, we can use
  the divisor equation to add a line condition, which introduces
  a factor $\frac{1}{d}$ and therefore leads to rational numbers.
  Then we can use both the tropical (see theorem \ref{thm-WDVV}) 
  and the classical WDVV equation (\cite{FP} or \cite{K}) and express
$\langle\tau_0(1)\prod_{k\in N}\tau_{r_k}(2) \rangle_d$  
 in terms of invariants 
$\langle\tau_0(0)^{l'}\tau_0(1)^{m'}\prod_{k\in N'}\tau_{r'_k}(2) \rangle_{d'}$  
 with $n' + \sum_{k \in N'} r'_k
  < n + \sum_{k \in N} r_k$. Repeating this procedure, we
  eventually end up with $n'=1$ (subcase I) or $r'_k = 0$ for all
  $k \in N'$, which is the 1. case.
\end{proof}

\section {Lattice paths} \label {latticepaths}

In this section, we present a lattice-paths algorithm to determine the numbers 
\begin{align*} \langle\prod_{k\in N} \tau_{r_k}(2)\rangle_d=
    \frac{1}{(d!)^3} \deg\left( \prod_{k \in N} \ev^*_k(P_k) \psi_k^{r_k}
      \cdot \M{n}{d}\right)
  \end{align*}
i.e.\ numbers of curves with Psi- and point conditions (and no line conditions; all other
numbers can be easily computed using string and divisor equation).
Note that in this case we need
\begin{displaymath}
 3d - 1 =  n + \sum_{k\in N} r_k
\end{displaymath}
to get a zero-dimensional cycle.

We use the fact that if we choose general point conditions, the intersection product
$\prod_{k \in N} \ev^*_k(P_k) \psi_k^{r_k}
      \cdot \M{n}{d}$ equals set-theoretically the set of all points corresponding to curves satisfying the Psi- and the point conditions (see \ref{rem-sandz}).
Each such curve $C$ has to be counted with weight, and it is counted with the weight $\frac{1}{d!^3}|\det_C(\ev)|$ (see theorem \ref{MultEqualDet}), where $\ev$ denotes the product of all evaluation maps (at both coordinates).
Note that no such curve can have a string since this would provide a deformation of the curve
described in the proof of \ref{ContractedEdge}, which contradicts \ref{GeneralConditions} (a).

We pick a certain configuration of points and count dual Newton subdivisions of curves passing through the points and satisfying the Psi-conditions. The dual Newton subdivisions are in fact dual to the image $h(\Gamma)\subset \RR^2$ of the graph in the plane. In particular, the labels of the non-contracted ends are lost.
That means we have to count tropical curves without labels on the non-contracted ends, and then multiply with the number of possibilities to set labels.

There is a map $c:\calM_{0,n}^{\lab}(\RR^2,d)\rightarrow \calM_{0,n}(\RR^2,d)$ which forgets the labels of the non-marked ends. This map is a cover, the number of preimages is the number of ways to set labels.
The biggest number of preimages is $d!^3$. However, not every point in a facet has this number of preimages:
the curve in $\calM_{0,0}(\RR^2,2)$ pictured below has only 4 preimages, not 8, since the two ends in direction $(-1,0)$ are not distinguishable.
\begin{center}
\begin{picture}(0,0)%
\includegraphics{notdist.pstex}%
\end{picture}%
\setlength{\unitlength}{3947sp}%
\begingroup\makeatletter\ifx\SetFigFont\undefined%
\gdef\SetFigFont#1#2#3#4#5{%
  \reset@font\fontsize{#1}{#2pt}%
  \fontfamily{#3}\fontseries{#4}\fontshape{#5}%
  \selectfont}%
\fi\endgroup%
\begin{picture}(4224,999)(4189,-3598)
\put(8026,-2836){\makebox(0,0)[lb]{\smash{{\SetFigFont{10}{13.2}{\familydefault}{\mddefault}{\updefault}{\color[rgb]{0,0,0}$\RR^2$}%
}}}}
\end{picture}%

\end{center}
Let $ C'\in \prod_{k \in N}  \psi_k^{r_k}
      \cdot \M{n}{d}$ and let $C$ be the curve after forgetting the labels of the non-contracted ends.
Assume that the facet $\sigma$ in $M_{0,n}(\RR^2,d)$ in which $C$ lies has $s$ preimages under the cover above.
Thus $C$ has to be counted with $\frac{s}{d!^3}|\det_{C'}(\ev)|$.
Assume that $C$ has $t$ vertices $V_1,\ldots,V_t$ such that $b_{ij}$ non-contracted ends of the same direction $-e_j$ are adjacent to $V_i$ (where $j$ goes from $0$ to $2$ and $e_0:=-e_1-e_2$). 
Then $\nu_C:=\prod_{i=1}^t \prod_{j=0}^2\frac{1}{b_{ij}!}=\frac{s}{d!^3} $ and we have to count $C$ with $\nu_C|\det_{C'}(\ev)|$.

First, we want to understand this weight locally in terms of vertex multiplicities. 
We define another weight that we denote by $\mult(C)$ and we show that it is equal to $\nu_C|\det_{C'}(\ev)|$.

\begin{definition}
Let $C'\in \prod_{k \in N}  \psi_k^{r_k}
      \cdot \M{n}{d}$ and let $C$ be the curve after forgetting the labels of the non-contracted ends.
Define the weight $\mult(C)$ as $\nu_{C}$ times the product of the multiplicities of those (necessarily 3-valent) vertices without any marked ends on them (see \cite{Mi03}, definition 2.16).
\end{definition}

\begin{example}\label{ex-mult}
Let $C$ be the curve as in the picture below. (For this example, we chose some other degree, not $d$, to keep the picture nice.) 
As in remark \ref{remark-basis} we choose coordinates to write down an explicit matrix for $\ev$.
Choose $V$ to be the root vertex. Then the matrix of $\ev$ is 
 \[ \left( \begin {array}{cccc}
           1 & 0 & v_{1,1} &  0  \\
           0 & 1 & v_{1,2} &  0  \\
           1 & 0 &  0  & v_{2,1} \\
           0 & 1 &  0  & v_{2,2}
         \end {array} \right). \]
The absolute value of the determinant is equal to $|\det(v_1,v_2)|$, which is the multiplicity of the $3$-valent vertex $V$.
\begin{center}
\input{multc.pstex_t}
\end{center}
\end{example}


\begin{lemma}\label{lem-mult}
Let $C'\in \prod_{k \in N}  \psi_k^{r_k}
      \cdot \M{n}{d}$ and let $C$ be the curve after forgetting the labels of the non-contracted ends.
Then $\nu_C|\det_{C'}(\ev)|=\mult(C)$ if $C$ has no string.
\end{lemma}

\begin{proof}
We have to show that $|\det_{C'}(\ev)|$ equals 
the product of multiplicities of all $3$-valent vertices without marked end.

This is an induction on the number of bounded edges. Curves with no bounded edge satisfy \linebreak $|\det_{C'}(\ev)| = \prod_V \mult(V)=1$ (the product is empty). Curves with $2$ bounded edges (as the one in the example \ref{ex-mult}) need to have the two marked ends at the two ``outer'' vertices, because otherwise there is a string. So as in the example, there is one ``interior'' 3-valent vertex without a marked end, and $|\det_{C'}(\ev)|$ is the multiplicity of the vertex.

Now we can assume we have 4 or more bounded edges. (The number of bounded edges is even, $2n-2$.) We choose one such that there are still bounded edges on both sides of it.
(If such an edge does not exist, it means we have a ``star-shaped'' tropical curve with one vertex in the middle and all bounded edges around. If one of those bounded edges was not adjacent to a marked end, $C$ has a string, so we can assume that all bounded edges are adjacent to a marked end. If there is no marked end in the middle, we then have $b+2=2b$ where $b$ is the number of bounded edges, since $b+2$ is the number of coordinates of this cone, and $2b$ is the number of coordinates of the $b$ marked ends. So $b=2$ and we are in the situation of the example. If there is a marked end in the middle, we have $b+2=2b+2$ so $b=0$, so this curve has no bounded edge and counts one. Now we can assume we do not have a star-shaped curve, and there is in fact a bounded edge with bounded edges on both sides.) 
 Then we cut this edge to get two curves $C_1$ and $C_2$.
Let $I_i$ denote the subset of marked ends on $C_i$ and let $e_i$ be the number of non-contracted ends of $C_i$. We make the cut edge a new non-contracted end of $C_i$, so $C_i$ has in fact $e_i+1$ non-contracted ends, one of them the special new end.
Assume \begin{displaymath}\# I_1\leq e_1-2-\sum_{k\in I_1}r_k,\end{displaymath} then if we remove all the closures of marked ends (as in lemma \ref{lem-string}) we get \begin{displaymath}\sum_{k\in I_1}r_k+\#I_1+1\end{displaymath}  connected components, which is less than or equal to \begin{displaymath}\sum_{k\in I_1}r_k+e_1-2-\sum_{k\in I_1}r_k+1=e_1-1.\end{displaymath}  So there must be a connected component which has two non-contracted ends of $C$ (not the new end of $C_1$). Hence $C$ has a string, which contradicts the assumption.
We have \begin{displaymath}\#I_1+\#I_2=n=3d-1-\sum r_k= e_1+e_2-1-\sum r_k.\end{displaymath} 
Therefore \begin{displaymath}\# I_1= e_1-\sum_{k\in I_1}r_k\mbox{ and }\# I_2= e_2-1-\sum_{k\in I_2}r_k\end{displaymath}  without restriction.
As in remark \ref{remark-basis}, we pick coordinates to write down an explicit matrix for $\ev$.
$C_1$ has \begin{displaymath}e_1+\#I_1-2-\sum_{k\in I_1}r_k=2\#I_1-2\end{displaymath}  bounded edges.
We pick the root vertex to be the boundary vertex of the cut edge in $C_1$. We order the basis elements such that the root vertex comes first, then the bounded edges in $C_{1}$, then the cut edge, then the bounded edges in $C_2$. We order the basis of $\RR^{2n}$ such that the marked ends in $C_1$ come first and then the marked ends in $C_2$.
Then the matrix of $\ev$ for $C$ is a block matrix. The block on the top left is just the matrix of $\ev$ for $C_1$ --- so by induction, the product of multiplicities of 3-valent unmarked vertices of $C_1$. The top right block is $0$, because no marked end on $C_1$ needs a bounded edge of $C_2$. 
The bottom right block has the same determinant as the matrix $\ev$ for $C_2$, when we add a marked end on the cut edge and make its end vertex the root vertex. So the determinant of this block is again by induction the product of multiplicities of 3-valent unmarked vertices in $C_2$. This proves the claim.
\end{proof}

Now we know that $\frac{1}{(d!)^3}\deg\left( \prod_{k \in N} \ev^*_k(P_k) \psi_k^{r_k}
      \cdot \M{n}{d}\right)$ equals the number of all curves $C\in \mathcal{M}_{0,n}(\RR^2,d)$ satisfying the Psi- and the point conditions, each counted with weight $\mult(C)$.
We want to simplify this count even further: we do not want to count parametrized tropical curves $C=(\Gamma,x_i,h)$, but we want to count their images in $\RR^2$. 

\begin{definition}
Let $C=(\Gamma,x_i,h)\subset \prod_{k \in N}  \psi_k^{r_k}
      \cdot \M{n}{d}$. In the image $h(\Gamma)$, some edges may lie on top of each other. Mark each edge in the image $h(\Gamma)$ by a partition reflecting the weights of all edges which map onto this image edge.
The image $h(\Gamma)$ together with those partitions is called the \emph{labelled image} of $C$.
\end{definition}
\begin{example}\label{ex-labelimage}
The following picture shows a tropical curve and its labelled image.
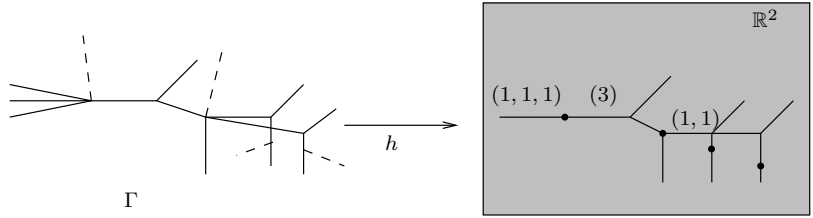
\begin{center}
\input{labelledimage.pstex_t}
\end{center}
\end{example}

Given a labelled image, there can be different possible parametrizations. Ambiguity may for example arise if the labelled image comes from a parametrization that maps vertices on top of each other. We could then also parametrize this labelled image with a graph where the two vertices are replaced by only one.
To avoid this ambiguity, we need a slightly more special notion of general conditions, which we call restricted general conditions. 
This definition is cooked up in such a way that we exactly avoid all ambiguity and make parametrizations unique.

\begin{definition}
A curve $C=(\Gamma,x_i,h)$ is called \emph{simple}, if it satisfies:
\begin{enumerate}
\item the map $h$ is injective on vertices,
\item if $h(V)\in h(e)$ for a vertex $V$ and an edge $e$ then $V$ is adjacent to an edge $e'$ which is mapped on top of $e$,
\item if two edges $e$ and $e'$ are mapped on top of each other, then they share a vertex,
\item assume $p\in \RR^2$ is a point through which more than two edges pass. Divide the edges into equivalence classes depending on the slope of the line to which they are mapped. Then we have at most $2$ equivalence classes.
\end{enumerate}
\end{definition}

\begin{definition}
The subset of $\RR^{2n}$ of \emph{restricted general conditions} is defined to be the subset of the set of general conditions such that only simple curves $C\in \prod_{k \in N}  \psi_k^{r_k}
      \cdot \M{n}{d}$ pass through the points (i.e.\ satisfy $\ev(C)=(P_1,\ldots,P_n)$).
\end{definition}
\begin{remark}
It is easy to see that the subset of restricted general conditions is still open and dense. Points which are not restricted general admit a non-simple curve. Being not simple sums up to codimension 1 conditions, hence only the image under $\ev$ of certain lower-dimensional subsets of $\prod_{k \in N}  \psi_k^{r_k}
      \cdot \M{n}{d}$ is not restricted general.
\end{remark}

\begin{lemma}\label{lem-labelledimageunique}
Given a labelled image of a tropical curve through restricted general conditions, there is exactly one abstract tropical curve $(\Gamma,x_i)$ and one map $h$ parametrizing this labelled image and sending the marked ends to the $P_i$.
\end{lemma}
\begin{proof}
Clearly there is a parametrization $C=(\Gamma,x_i,h)$ of the labelled image, we just need to show that it is unique.
Since the $P_i$ are general, $C$ cannot have a contracted bounded edge. 
If it had a contracted bounded edge, we could vary the length of this edge without changing the image, in contradiction to \ref{GeneralConditions}(a).
Hence all edges can be seen in the image $h(\Gamma)$.
Due to the conditions for being simple, we can also distinguish the images of vertices and the images of all edges in the labelled image. Because of the labels we know whether edges lie on top of each other. If there are edges lying on top of each other, then we know that they have to share a vertex.
If there is a vertex $V$ with two edges of the same direction, it has to be more than 3-valent.
If it was $3$-valent, then by the balancing condition the 3 edges would be mapped to a line. At least one of the 3 edges is bounded, and so we could change the length of this edge (and accordingly the lengths of the other two edges, if necessary) without changing the image of the curve. That contradicts \ref{GeneralConditions}(a).
Since $V$ is more than $3$-valent, there must be a marked end adjacent to it. It is not possible that $2$ (or more) of the points $P_i$ lie on a line with a direction that can be the direction of an edge. Thus if we have two edges in the labelled image on top of each other, there must be exactly one adjacent vertex which passes through a point $P_i$. Thus we know that the edges have to be connected at that vertex when we built the parametrization. 
\end{proof}

\begin{definition}
Let $C$ be a curve in $\prod_{k \in N}  \psi_k^{r_k}
      \cdot \M{n}{d}$ passing through restricted general conditions.
Draw a dual Newton subdivision to the image $h(\Gamma)$ and label the dual edges also with the partitions belonging to the edges of the labelled image $h(\Gamma)$.
This is called a \emph{labelled dual Newton subdivision}.
Mark the polygons dual to vertices which are adjacent to a marked end $x_i$.
Those marked polygons in $\Delta_d$ together with the partitions belonging to their boundary edges is called the \emph{set of dual marked polygons of $C$}. 
\end{definition}

\begin{example}
The following picture shows the labelled dual Newton subdivision to the labelled image from example \ref{ex-labelimage}.
Next to it, we can see the set of dual marked polygons of $C$.
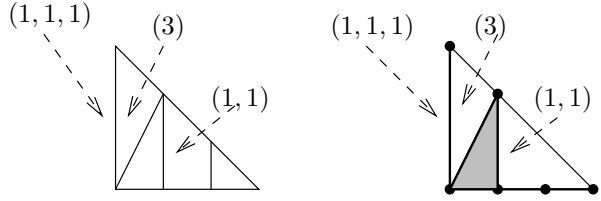
\begin{center}
\input{exlabeldual.pstex_t}
\end{center}
\end{example}

Our aim is to count dual marked polygons to curves in $\prod_{k \in N}  \psi_k^{r_k}
      \cdot \M{n}{d}$. To do that, we have to choose a special point configuration $\mathcal{P}$ as our condition. This configuration is chosen is such a way that the set of dual marked polygons can be described as something like a generalized lattice path that we call a \emph{rag rug}.
We will now first introduce labelled lattice paths and rag rugs, and then show that the count of rag rugs equals the count of labelled images of curves in $\prod_{k \in N}  \psi_k^{r_k}
      \cdot \M{n}{d}$ passing through $\mathcal{P}$ (with weight $\mult(C)$).

 Let $\Delta_d$ be the triangle with endpoints $(0,0)$, $(d,0)$ and $(0,d)$. 
Fix $\lambda$ to be a linear map of the form \begin{displaymath} \lambda:\RR^2 \to\RR:\;(x,y) \mapsto x-\varepsilon y ,\end{displaymath} where $ \varepsilon $ is a small irrational
number. 
Recall that a path $\gamma: [0,n] \rightarrow\RR^2$ is called a lattice path if
  $\gamma |_{[j-1,j]}, j=1,\ldots,n$ is an affine-linear map and $\gamma(j)
  \in\ZZ^{2} $ for all $ j=0\ldots,n $.
For $n=1,\ldots,n$, we call $\gamma |_{[j-1,j]}([j-1,j])$ a \emph{step} (the $j$-th step) of the lattice path $\gamma$.
A lattice path is called $\lambda$-increasing if $\lambda \circ \gamma$
is strictly increasing. Let $p:=(0,d) $ and $q:=(d,0) $ be the points in $
\Delta:=\Delta_d $ where $\lambda|_{\Delta}$ reaches its minimum (resp.\
maximum). 
Let $G$ be a line in $\RR^2$ orthogonal to $\ker(\lambda)$. Then $G$ divides the plane into two halfplanes. We will denote the upper one by $H^+$ and the lower one by $H^-$.

\begin{definition}
A \emph {labelled $\lambda$-increasing lattice path} in $\Delta$ is a $\lambda$-increasing lattice path from $p$ to $q$ such that the $k$-th step is labelled by a partition $\alpha_k=((\alpha_k)_1,\ldots,(\alpha_k)_{r_k})$ of the integer length of this step, that is \ $(\alpha_k)_1+\ldots+(\alpha_k)_{r_k}=\#(\ZZ^2\cap \gamma([k-1,k]))-1 $.
\end{definition}

\begin{remark}\label{rem-deltas}
Let $\delta$ be a labelled $\lambda$-increasing lattice path from $p$ to $q$ whose image is contained in the boundary $\partial \Delta$ and whose steps are labelled with partitions consisting of only ones. All those paths will be possible end paths for the recursion defining multiplicity.
The following picture shows 3 examples for $\Delta_3$.
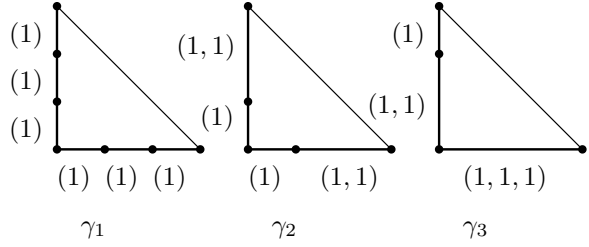
\begin{center}
\input{endpaths.pstex_t}
\end{center}

\end{remark}

\begin{definition}\label{def-mult}
We define the \emph{positive multiplicity} $\mu_+$ (resp.\ \emph{negative multiplicity} $\mu_-$) of a labelled $\lambda$-increasing lattice path recursively as follows:
\begin{enumerate}
\item For a possible end path $\delta$ as in remark \ref{rem-deltas} going \emph{clockwise} from $p$ to $q$ (resp.\ \emph{counterclockwise}) with $n$ steps we define $\mu_\pm(\delta):=\prod_{k=1}^n 1/(|\alpha_k|!)$, where $|\alpha_k|$ denotes the size of the partition of the $k$-th step (recall it has to be a partition with only ones as entries).

\item For a labelled $\lambda$-increasing lattice path $\gamma$ which is not a possible end path, assume that the $k$-th and the $k+1$-th step form the first \emph{left} (resp.\ \emph{right}) corner of the path $\gamma$. (If no such turn exists, we define $\mu_\pm(\gamma):=0$.)
\\
\\
Define a finite set of lattice paths as follows:
\begin{itemize}
\item pick an integer $r$ with $0\leq r \leq \min\{| \alpha_k|,| \alpha_{k+1}| \}$,

\item pick a set $S$ of $r$ pairs \begin{displaymath}S=\big\{[(\alpha_k)_{i_1},(\alpha_{k+1})_{j_1}],\ldots,[(\alpha_k)_{i_r},(\alpha_{k+1})_{j_r}]\big\}\end{displaymath}
such that the multiset $\{(\alpha_k)_{i_1},\ldots,(\alpha_k)_{i_r}\}$ 
is a subset of the multiset $\{(\alpha_k)_{1},\ldots,(\alpha_k)_{r_k}\}$ 
and the multiset $\{(\alpha_{k+1})_{j_1},\ldots,(\alpha_{k+j})_{j_r}\}$ 
is a subset of the multiset \linebreak
$\{(\alpha_{k+1})_{1},\ldots,(\alpha_{k+1})_{r_{k+1}}\}$. 
\end{itemize}

For each $l=1,\ldots,r$, build a triangle $T_{r,S,l}$ with one edge of integer length $(\alpha_k)_{i_l}$ and one edge of integer length $(\alpha_{k+1})_{j_l}$ (in the direction of the $k$-th resp.\ $k+1$-th step).
Let $M_{r,S}$ be the Minkowski sum of all triangles $T_{r,S,l}$ for $l=1,\ldots,r$, and edges $e_s$ in direction of the $k$-th step of integer length $(\alpha_k)_{s}$ for all $s$ which are not one of the $i_l$ and edges $f_t$ in direction of the $k+1$-th step of integer length $(\alpha_{k+1})_{t}$ for all $t$ which are not one of the $j_l$. Label each edge $E$ of $M_{r,S}$ with a partition reflecting the integer lengths of edges $e_s$, $f_t$, and edges of triangles $T_{r,S,l}$ that contribute to $E$.
Think of the polygon $M_{r,S}$ as sitting in the corner built by step $k$ and $k+1$ of $\gamma$, and define a new labelled $\lambda$-increasing lattice path ${\gamma}_{r,S}$ by going the other way around $M_{r,S}$.
If $M_{r,S}$ does not fit inside the polygon $\Delta$, we define $\mu_{\pm}(\gamma_{r,S})=0$.
The positive multiplicity of this new labelled $\lambda$-increasing lattice path is known recursively, because it includes a smaller area with the possible end paths.
We define
\begin{displaymath}
\mu_\pm(\gamma)=\sum_r \sum_{S}  \Area(T_{r,S,1})\cdot \ldots\cdot \Area(T_{r,S,r})\cdot \mu_\pm(\gamma_{r,S}),
\end{displaymath}
(where $\Area(T)$ is the normalized lattice area, i.e.\ the area of the simplex with vertices $(0,0)$, $(1,0)$ and $(0,1)$ is defined to be $1$).
\end{enumerate}

\end{definition}

\begin{example}
For the $3$ possible end paths from remark \ref{rem-deltas}, we have multiplicity $\mu_-(\gamma_1)=1$, $\mu_-(\gamma_2)=\frac{1}{4}$ and $\mu_-(\gamma_3)=\frac{1}{12}$.
\end{example}

\begin{example}
The following picture shows an example of the recursion from definition \ref{def-mult} to compute the positive multiplicity of a labelled path $\gamma$ in $\Delta_3$. The first left turn is from step $2$ to step $3$. 
We have $3$ choices for $r$: $r=0$, $r=1$ or $r=2$, since both the partition of step $2$ as the partition of step $3$ contain $2$ elements. No matter what we choose for $r$, there is just one choice for the set $S$ ($r$ pairs consisting of all ones), since both partitions contain only ones.

For $r=0$ and $S=\emptyset$, $M_{0,\emptyset}$ is a square of size $2$ which does not fit inside ${\Delta}_3$. Therefore the multiplicity of $\gamma_{0,\emptyset}=0$.

For $r=1$ and $S=\{(1,1)\}$, $M_{1,S}$ is a pentagon. The integer length of each new side is one. The new side of direction $(0,1)$ is labelled by the partition $(1)$, because it comes from the edge $e_s$ in direction of the $2$-nd step of integer length $1$ which is not one of the $i_l$ in the set of pairs $S$. The new side of $M_{1,S}$ of direction $(1,-1)$ comes with label $(1)$, because it comes from a side of the triangle $T_{1,S,1}$ of integer length $1$. The side of direction $(1,0)$ comes from an edge of the $3$-rd step which is not part of $S$, and gets label $(1)$ as well. The area of $T_{1,S,1}$ is one.

For $r=2$, we have $S=\{(1,1),(1,1)\}$, and $M_{2,S}$ is a triangle of size $2$ whose new side gets the label $(1,1)$ because it comes from $2$ sides of the two triangles $T_{2,S,1}$ and $T_{2,S,2}$.
The area of both triangles $T_{2,S,1}$ and $T_{2,S,2}$ is one.

The picture shows how the recursion goes on after the first step. The choices where $r=0$ for ${\gamma}_{1,S}$ or where $r=1$ for $\gamma_{2,S}$ are left out because they yield to a path of multiplicity $0$. We end up with one path of multiplicity $1$ and one of multiplicity $\frac{1}{2}$, so $\mu_+(\gamma)=\frac{3}{2}$.

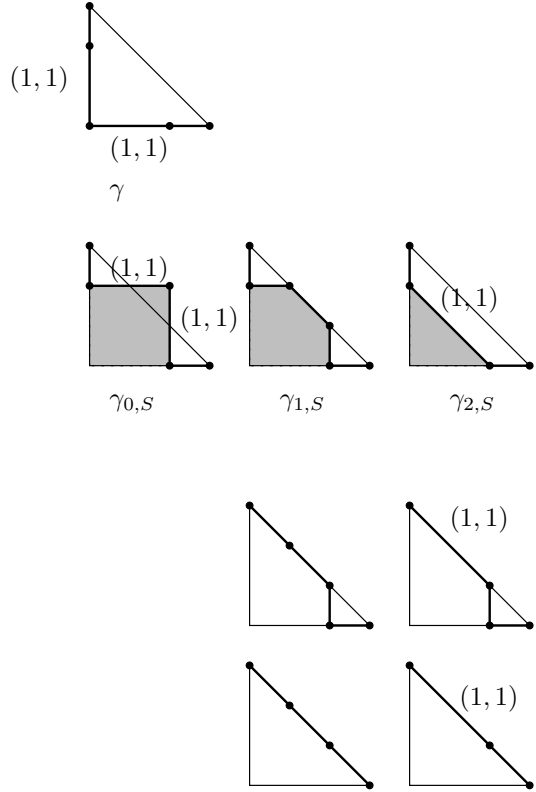
\begin{center}
\input{multpsiweg.pstex_t}
\end{center}
\end{example}

\begin{definition}\label{def-hm}
Let $F$ be a set of $n$ convex polytopes $Q_1,\ldots,Q_n$ inside $\Delta$ whose endpoints are lattice points of $\Delta$ and whose boundary edges $e$ are labelled by partitions. It is possible that a polygon $Q_i$ is $1$-dimensional, i.e.\ just an edge itself, then it has two partitions as labels, one for each outward pointing normal vector.
We call $F$ a \emph{rag rug} of the form $(r_1,\ldots,r_n)$ if the following conditions are satisfied:
\begin{enumerate}
\item the (outside) label $\alpha_e$ of an edge $e$ in the boundary of $\Delta_d$ is $\alpha_e=(1,1,\ldots,1)$,
\item two polygons $Q_i$ and $Q_j$ intersect in at most one point,
\item boundary edges whose outward normal vector points into $H^+$ (starting at $G$) (with their corresponding labels) form a labelled $\lambda$-increasing lattice path from $p$ to $q$ that we will denote by $\gamma^+$,
\item boundary edges whose outward normal vector points into $H^-$ (with their corresponding labels) form a labelled $\lambda$-increasing lattice path from $p$ to $q$ that we will denote by $\gamma^-$,
\item the order of the polytopes $Q_1,\ldots,Q_n$ agrees with the obvious order given by the paths $\gamma^+$ resp.\ $\gamma^-$,
\item the sum of the sizes of the partitions of the boundary edges of $Q_i$ is equal to $r_i+2$,\begin{displaymath}
\sum_{e\mid e\mbox{ edge of }Q_i}| \alpha_e|=r_i+2.
                                                                                              \end{displaymath}
\end{enumerate}
We define the multiplicity $\mu(F)$ to be $\mu_{+}(\gamma_+)\cdot \mu_-(\gamma_-)$.
\end{definition}

\begin{example}
The following picture shows a rag rug $F$ of the form $(2,2,0,0)$ in $\Delta_3$, and the two labelled $\lambda$-increasing lattice paths $\gamma_+$ and $\gamma_-$. 
For all edges of integer length one, the corresponding partitions are just $(1)$ and we did not mark this in the picture.
We have $\mu_+(\gamma_+)=3$ and $\mu_-(\gamma_-)=\frac{1}{6}$, so $\mu(F)=\frac{1}{2}$.
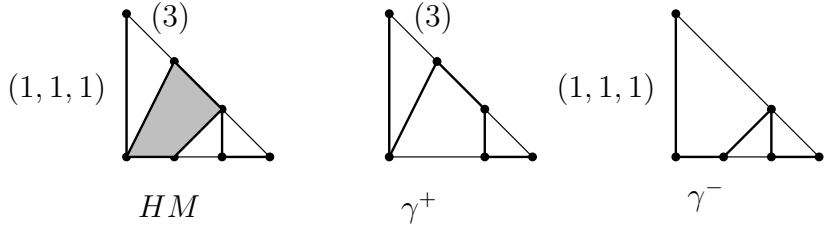
\begin{center}
\input{psiwege1.pstex_t}
\end{center}
\end{example}

\begin{definition}
Given $d$, $n$ and numbers $(r_1,\ldots,r_n)$ we define $N_{\rr}(d,n,(r_1,\ldots,r_n))$ to be the number of rag rugs of form $(r_1,\ldots,r_n)$, counted with multiplicity as defined in \ref{def-hm}.
\end{definition}

\begin{remark}
Note that this definition generalizes Mikhalkin's lattice path count (see \cite{Mi03}). A $\lambda$-increasing lattice path $\gamma$ from $p$ to $q$ is a rag rug of form $(0,\ldots,0)$. We have to attach labels $(1)$ to each edge. 
The two paths $\gamma_+$ and $\gamma_-$ agree with $\gamma$.
In the recursion for the lattice path count, we define $\mult_\pm(\gamma)$ depending on the multiplicity of two other paths $\gamma'$ and $\gamma''$. $\gamma'$ is the path that cuts the corner, and $\gamma''$ is the path that completes the corner to a parallelogram.
In our definition, we can choose $r=0$ or $r=1$. For $r=0$, we have $S=\emptyset$ as only choice. The polygon $M_{0,\emptyset}$ is the parallelogram which is equal to the Minkowski sum of the two steps of the corner. For $r=1$, we have $S=\{(1,1)\}$ as only choice, and $M_{1,S}$ is the triangle formed by the two steps of the corner.
Since all partitions are just $(1)$, also the end paths have only those partitions, so that there is in fact only one end path, the path $\delta_\pm$. It has multiplicity $1$.
Therefore our definition gives the same multiplicity in this case.
\end{remark}

It is not true that $N_{\rr}(d,n,(r_1,\ldots,r_n))= \frac{1}{(d!)^3}\deg\left( \prod_{k} \ev^*_k(P_k) \psi_k^{r_k}
      \cdot \M{n}{d}\right)$, since we count also reducible curves with the rag rugs. 
     
We therefore have to modify the count and throw away the dual subdivisions corresponding to reducible tropical curves.

\begin{definition}
Given a rag rug $\gamma$ and the two corresponding lattice paths $\gamma_+$ and $\gamma_-$, perform the recursion to compute their multiplicity and keep track of the polygons $M_{r,S}$ that the new paths $\gamma_{r,S}$ in the recursion enclose with $\gamma_\pm$. This way we end up with a set of labelled Newton subdivisions. We call this the set of \emph{possible labelled Newton subdivisions for $\gamma$}. 
The recursion allows us to assign a multiplicity to a possible labelled Newton subdivision, so that the multiplicity of $\gamma$ is equal to the sum of the multiplicities of the possible labelled Newton subdivisions for $\gamma$.
\end{definition}
\begin {definition}\label{def-multNewtonsub}
Given a labelled Newton subdivision, draw a dual labelled image and then the unique tropical curve mapping to this image. This is well-defined up its position in $\RR^2$ and the lengths of its bounded edges.
We say that the Newton subdivision is \emph{reducible} if the tropical curve mapping to a dual labelled image is \emph{reducible} (again, this does not depend on the choice of dual labelled image).
Otherwise, we say it is irreducible.
\end {definition}

\begin{remark}
It is possible to express the reducibility condition in terms of the Newton subdivision itself and not in terms of the dual tropical curve.
A labelled marked Newton subdivision is reducible if and only if it admits a mixed subdivision where the marked polygons $Q_i$ remain unmixed (i.e.\ come from sums of the form $Q_i+v_1+\ldots+v_t$, where the $v_j$ are vertices of the subdivision of the $j$-th summand). For details, see \cite{M08}.
\end{remark}

\begin{definition}\label{def-irrHM}
For a rag rug $\gamma$ as defined in \ref{def-hm} define its irreducible multiplicity $\mult'(\gamma)$ to be the multiplicity $\mult(\gamma)$ minus the number of possible reducible Newton subdivision (counted with multiplicity).
Given $d$, $n$ and numbers $(r_1,\ldots,r_n)$ we say $N'_{\rr}(d,n,(r_1,\ldots,r_n))$ is the number of rag rugs counted with their irreducible multiplicity.
\end{definition}

\begin{definition}
Let $\mathcal{P}_\lambda=(P_1,\ldots,P_n)$ denote $n$ restricted general point conditions on the line $G$ orthogonal to $\ker(\lambda)$ such that the distance between $P_i$ and $P_{i+1}$ is much bigger than the distance between $P_{i-1}$ and $P_i$.
\end{definition}
\begin{lemma}\label{lem-dualHM}
Let $C\in \prod_{k} \psi_k^{r_k}
      \cdot \M{n}{d}$ with $\ev(C)=\mathcal{P}_\lambda$. Then the set of dual marked polygons of $C$ is a rag rug of the form $(r_1,\ldots,r_n)$.
\end{lemma}
\begin{proof}
The polygon $Q_i$ dual to $P_i$ is convex and has to satisfy \begin{displaymath}\sum_{e\mid e\mbox{ edge of }Q_i}| \alpha_e|=r_i+2,\end{displaymath} where $\alpha_e$ denotes the partition belonging to $e$. This is true since the marked end $x_i\subset \Gamma$ is adjacent to a vertex of valence $r_i+3$ and we can see all edges (except the contracted end $x_i$) in the labelled image (and thus in their labelled dual Newton subdivision, too).
(Outside) labels of edges in the boundary of $\Delta_d$ have only ones as entries, since the ends of $C$ are all of weight $1$.
That the boundaries of those polygons form labelled $\lambda$-increasing lattice paths follows analogously to \cite{Mi03}, 8.27 (or \cite{Ma06}, 5.48 for more details). 
\end{proof}

\begin{theorem}
The number $N'_{\rr}(d,n,(r_1,\ldots,r_n))$ from definition \ref{def-irrHM} equals the intersection product $\frac{1}{(d!)^3}\deg\left( \prod_{k} \ev^*_k(P_k) \psi_k^{r_k}
      \cdot \M{n}{d}\right)=\langle \prod_k \tau_{r_k}(2) \rangle$.
\end{theorem}

\begin{proof}
To determine $\frac{1}{(d!)^3}\deg\left( \prod_{k} \ev^*_k(P_k) \psi_k^{r_k}
      \cdot \M{n}{d}\right)$, we can draw all labelled images of tropical curves that pass through $\mathcal{P}_\lambda$ and count them each with their weight $\mult(C)$ which is $\nu_C$ times the product of the multiplicities of non-marked vertices. We show that this count is equivalent to counting irreducible possible labelled Newton subdivision for all rag rugs of the form $(r_1,\ldots,r_n) $.

The proof is a generalization of the proof of theorem 2 of \cite{Mi03}.
We know that each $C$ leads to a rag rug $\gamma_C$ as in lemma \ref{lem-dualHM}.
Each rag rug yields a set of possible labelled Newton subdivision, with multiplicity.
We will show that for each such possible labelled Newton subdivision, there is a dual tropical curve $C$ through $\mathcal{P}_\lambda$ of the same weight.
At the same time, we show that for each curve $C$ through $\mathcal{P}_\lambda$, the dual labelled Newton subdivision is possible for the rag rug $\gamma_C$.

Let $\gamma$ be a rag rug. 
The recursion for $\gamma_+$ yields possible subdivision of $\Delta_d$ above the polygons $Q_i$. They correspond to the part of a tropical curve above $G$. Analogously, possible subdivisions for $\gamma_-$ correspond to the parts of tropical curves below $G$. 
From lemma \ref{lem-mult} it follows that weight of a tropical curve can be computed locally, so the weight of $C$ is equal to the weight of the part above $G$ times the weight of the part below $G$. The same is true for the multiplicity of the dual Newton subdivision.
Therefore it is enough to show that for each subdivision above the $Q_i$, there is a dual part of a tropical curve above $G$ of the same weight, and that each part of a tropical curve above $G$ is dual to a possible subdivision above the $Q_i$. The corresponding statement for subdivisions below the $Q_i$ and parts of tropical curves below $G$ follows analogously, and thus the complete statement follows.

For each point $P_i\in \mathcal{P}_\lambda$, draw edges emanating from $P_i$ of directions dual to the boundary edges of $Q_i$ and with the same partitions as labels. Draw a line $G'$ in $H_+$ parallel to $G$, such that the strip between $G$ and $G'$ encloses one intersection of the edges we have drawn through the $P_i$. This intersection of edges corresponds to the first left turn of the path $\gamma_+$, since the distances between the $P_i$ are increasing. 
Let us determine the possibilities how the tropical curve can go on at this point. We should think about both edges as a set of edges of weights given by the partition. Edges can either meet in a $3$-valent vertex, or intersect. First, we pick $r$ less than the smaller number of edges in a set to determine how many edges should meet in a $3$-valent vertex. Then we pick a set of $r$ pairs of weights to determine which edges should meet in a $3$-valent vertex. The other edges intersect. 
The weight locally in the strip between $G$ and $G'$ is equal to the product of areas of triangles dual to the $3$-valent vertices because of lemma \ref{lem-mult}. The dual polygon is the Minkowski sum of those triangles and the remaining edges which intersect.
Therefore the recursion for the multiplicity of $\gamma_+$ corresponds to the possibilities for a labelled image of a tropical curve in the strip between $G$ and $G'$ and keeps track of the weight.
The end paths which do not have zero multiplicity are exactly those dual to ends of direction $(1,1)$ and weight one. The multiplicity of such an end path corresponds to the correction factor $\nu_C$ with which we have to divide the weight of a tropical curve if more than one non-contracted end is adjacent to the same vertex.
\end{proof}

\begin{example}
The following picture shows how to count $\langle \tau_2(2)^2\tau_0(2)^2\rangle_3$ using rag rugs. The left column shows all rag rugs of the form $(2,0,2,0)$ in the triangle $\Delta_3$. The middle column shows the possible Newton subdivisions for the rag rugs and their multiplicity. The third column shows sketches of the dual tropical curves.
\begin{center}
\input{zaehlbsp1.pstex_t}
\end{center}

\begin{center}
\input{zaehlbsp22.pstex_t}
\end{center}
\begin{center}
\input{zaehlbsp23.pstex_t}
\end{center}
\end{example}

\end {document}